\documentclass[review,11pt]{elsarticle}
\usepackage[utf8]{inputenc}
\usepackage[hidelinks]{hyperref}

\usepackage{microtype}
\usepackage[margin=2.5cm]{geometry}
\usepackage{setspace}
\usepackage[english]{babel}
\usepackage{color}
\usepackage{booktabs}
\usepackage{longtable}
\usepackage{graphicx}

\usepackage[ruled,vlined, linesnumbered]{algorithm2e}
\usepackage{algpseudocode}
\usepackage{empheq}
\usepackage[most]{tcolorbox}
\usepackage{ragged2e}

\usepackage[font=scriptsize,justification=justified,format=plain]{caption}
\usepackage{enumitem}
\usepackage{amsmath}
\usepackage{amssymb}
\usepackage{float}

\newtheorem{proposition}{Proposition}
\newtheorem{lemma}{Lemma}
\newtheorem{definition}{Definition}

\usepackage[normalem]{ulem}
\useunder{\uline}{\ul}{}
                    
\newcommand{\Zebra}{\texttt{Zebra}}
\newcommand{\Popstar}{\texttt{PopStar}}
\newcommand{\BandC}{\textit{branch-and-cut}}
\newcommand{\Avella}{\texttt{AvellaB\&C}}
\newcommand{\of}{\mathcal S}

\journal{}
\biboptions{authoryear}

\begin{document}
    \justifying
    \begin{frontmatter}
    
        \title{An efficient Benders decomposition for the \texorpdfstring{$p$}{p}-median problem 
        }

        \author[mymainaddress,mysecondaryaddress,mythirdaddress]{Cristian Durán\texorpdfstring{\corref{mycorrespondingauthor}}{}}
        \texorpdfstring{\ead{cristian.duran@ensta-paris.fr}}{}
        
        \author[mymainaddress,mysecondaryaddress]{Zacharie Alès}
        \ead{zacharie.ales@ensta-paris.fr}
        
        \author[mymainaddress,mysecondaryaddress]{Sourour Elloumi}
        \ead{sourour.elloumi@ensta-paris.fr}
        
        \address[mymainaddress]{UMA, ENSTA Paris, Institut Polytechnique de Paris, 828 Boulevard des Maréchaux, 91120~Palaiseau, France.}
        \address[mysecondaryaddress]{CEDRIC,  Conservatoire National des Arts et Métiers, 292 rue Saint-Martin, 75003~Paris, France.}
        \address[mythirdaddress]{LDSPS, Industrial Engineering Department, University of Santiago of Chile,  Avenida Víctor Jara 3769, 9160000~Santiago, Chile.
        }

        \begin{abstract}
            The $p$-median problem is a classic discrete location problem with several applications. It aims to open $p$ sites while minimizing the sum of the distances of each client to its nearest open site. We study a  Benders decomposition of the most efficient formulation in the literature. We prove that the Benders cuts can be separated by a polynomial time algorithm. The Benders decomposition also leads to a new compact formulation for the $p$-median problem. We implement a \emph{branch-and-Benders-cut} approach that outperforms state-of-the-art methods on benchmark instances by an order of magnitude.
        \end{abstract}
    

    \end{frontmatter}

    \section{Introduction} 
    
        Discrete location problems aim at choosing a subset of locations, from a finite set of candidates, in which to establish facilities in order to serve a finite set of clients. The sum of the fixed costs of the facilities and the allocation costs of supplying the clients must be minimized. 
        
        The $p$-median problem $(pMP)$ is an important location problem in which a given number $p$ of locations, usually called medians, have to be chosen from the set of candidate sites. In the $(pMP)$, no fixed costs are considered and the allocation costs are equal to the distance between clients and sites. More formally, given a set of $N$ clients $\{C_1,...,C_N\}$ and a set of $M$ potential sites to open $\{F_1,...,F_M\}$, let $d_{i j}$ be the distance between client $C_i$ and site $F_j$ and $p \in \mathbb{N}$ the number of sites to open. The objective is to find a set $S$ of $p$ sites such that the sum of the distances between each client and its closest site in $S$ is minimized. The $(pMP)$ leads to applications where the sites correspond to warehouses, plants, shelters, etc.  More recent applications can also be found in clustering processes in databases, where sub-groups of objects, variables, persons, etc. are identified according to defined criteria. We refer to~\cite{Marin2019} for a review on applications and resolution methods of the $(pMP)$.
        
        A great interest in solving large location problems has led to the development of various heuristics and meta-heuristics in the literature. However, the exact resolution of large instances remains a challenge. Some location problems have recently been efficiently solved using the Benders decomposition method within a \BandC\;approach (see e.g.,~\cite{Fischetti2017,Cordeau2019benders,gaar2021scaleable}).
        
        \subsection{Contribution and outline}
         
            In this paper, we explore a Benders decomposition for the $(pMP)$. We propose a polynomial time separation algorithm of the Benders cuts for the $(pMP)$. As a byproduct, we also obtain a new compact formulation of the $(pMP)$. We implement this decomposition together with other improvements within a \textit{branch-and-Benders-cut} approach. We show that our approach provides better results than the best exact resolution algorithm in the literature \Zebra~(\cite{Garcia2011}). We present our results on about 200 benchmark instances of different sizes (up to 115475 clients and sites) satisfying or not the triangle inequalities.
                
            The rest of the paper is organized as follows. Section~\ref{sec:literature} presents the literature review of $(pMP)$. Section~\ref{sec:Benders} describes our Benders decomposition method. Section~\ref{sec:exper} presents the computational results. In Section~\ref{sec:conclusions} we draw some conclusions together with research perspectives.
        
    \section{Literature review}\label{sec:literature}
       
        The $(pMP)$ was introduced by~\cite{Hakimi} where the problem was defined on a graph such that a client can only be assigned to an open neighbor site. It was showed that $(pMP)$ is an NP-hard problem by~\cite{Kariv1079}. The following is a summary of  the main formulations of this problem and its state-of-the-art exact resolution methods.
        
            \subsection{MILP formulations}\label{subsec:formulations}

            The classical mathematical programming formulation for the $(pMP)$ was proposed by \cite{revelle} who formulated the problem with a binary variable $y_j$ for each sites $F_j$  equals to 1 if the site is open and 0 otherwise; and a binary variable $x_{ij}$ equals to 1 if client $C_i$ is assigned to site $F_j$ and 0 otherwise. In the following, we denote by $[n]$ the set $\{1, 2, \hdots, n\}$ for any $n\in\mathbb N^*$.
            
            \begin{align}
                     &&\displaystyle \min \; \sum\limits_{i=1}^{N} \sum\limits_{j=1}^{M} d_{ij}x_{ij}    && \label{eq:f11}\\
               (F1): && \textrm{s.t.} \hspace{1.2cm} \sum\limits_{j=1}^{M} y_j     & =  p           &    &  \label{eq:f12} \\
                     && \sum\limits_{j=1}^M x_{ij}                                 &=  1            &  i &\in [N]   \label{eq:f13}\\
                     &&  x_{ij}                                                    & \leq y_j       &  i &\in [N], \ j \in [M] \label{eq:f14} \\
                     &&  x_{ij}                                                    & \geq 0         &  i &\in [N], \ j \in [M] \label{eq:f15} \\
                     \nonumber&&  y_{j}                                                     & \in \{ 0, 1\}  &  j &\in [M] 
            \end{align}
            Constraint~(\ref{eq:f12}) fixes the number of open sites to $p$. Constraints (\ref{eq:f13}) ensure that each client is assigned to exactly one site and Constraints (\ref{eq:f14}) ensure that no client is assigned to a closed site. The binary variables $x_{ij}$ can actually be relaxed as in Constraints~\eqref{eq:f15}. 
            
            An alternative formulation was proposed by~\cite{Cornuejols1980} which order for each client all its distinct distances to sites. More formally, for any client $i \in [N]$, let $K_i$ be the number of different distances from $i$ to any site. It follows that $K_i \leq M$. Let $D_i^1 < D_i^2 <...< D_i^{Ki}$ be these distances sorted. Formulation $(F2)$ uses the same $y$ variables as in formulation $(F1)$ and introduces new binary variables $z$. For any client $i \in [N]$ and $k \in [K_i]$, $z_{i}^k = 0$  if and only if there is an open site at distance at most $D_i^{k}$ from client $i$.
               
            \begin{align}
                && \displaystyle \min  \sum\limits_{i=1}^{N} \Big(  D_i^1  + \sum\limits_{k=1}^{K_i-1} (D_i^{k+1}- & D_i^{k})z_i^k \Big)  && \label{eq:f21}\\
                (F2):&& \textrm{s.t.} \hspace{2.8cm}   \sum\limits_{j=1}^{M} y_j   & =  p    \label{eq:f22} && \\
                && z_i^k +  \sum\limits_{j:d_{ij} \leq D_i^k} y_{j} &  \geq   1   &  i &\in [N], \ k \in [K_i]  \label{ct:zikd}\\
                && z_i^{k}        &  \geq   0      &  i & \in [N],\ k \in [K_i]  \label{ct:zik_positive}\\
                \nonumber     && y_{j}          & \in \{ 0, 1\}  &  j & \in [M] 
            \end{align} 
             \vspace{0.2cm}
            
             Objective~\eqref{eq:f21} minimizes the sum of the assignment distances over all clients. Constraints~\eqref{ct:zikd}  ensure that variable $z_i^k$ takes the value 1 if there is no site at a distance less than or equal to $D_i^k$ of client $i$. In that case  $(D_i^{k+1}-D_i^{k})$ is added to the objective. Otherwise, given the positive coefficients in the objective function, $z_i^k$ takes the value 0. Here again, the binary variables $z_i^k$ can be relaxed as in Constraints~\eqref{ct:zik_positive}.
             
             Formulation $(F2)$ can be much smaller than $(F1)$ and both have the same linear relaxation value~(\cite{Cornuejols1980}). In $(F1)$ the number of $x$ variables  is $N \times M$ and the number of constraints is $1 + N\times(1 + M )$. In $(F2)$, the number of $z$ variables is equal to $K =\sum\limits_{i=1}^{N} K_i$ and the number of constraints is equal to $1+K$. As $K \leq N \times M$, it follows that $(F2)$ has at most as many variables as $(F1)$ and at least $N$ less constraints. If the constraints coefficients matrix is sparse or if many facilities are equidistant from a given client then $K$ can be significantly smaller than $N\times M$.
        
            \cite{Elloumi2010} introduced an improved formulation based on $(F2)$. Given that, by definition,  $z_i^{k-1}$ equal to~0 implies that $z_i^{k}$ is also equal to~$0$, Constraints (\ref{ct:zikd})  can be replaced by (\ref{eq:f33}) and (\ref{eq:f34}). 
            \begin{align}
                &     &  \displaystyle \min  \sum\limits_{i=1}^{N} \Big(  D_i^1  + \sum\limits_{k=1}^{K_i-1} &(D_i^{k+1}- D_i^{k})z_i^k \Big)  && \label{eq:f31}\\
                &     & \textrm{s.t.} \hspace{1.5cm} \sum\limits_{j=1}^{M} y_j   & =  p   && \label{eq:f32}  \\
                         &(F3):& z_i^1 +  \sum\limits_{j:d_{ij} = D_i^1} y_{j} &  \geq    1   &  i &\in [N]  \label{eq:f33} \\ 
                         &     & z_i^k +  \sum\limits_{j:d_{ij} = D_i^k} y_{j} &  \geq   z_i^{k-1}    &  i &\in [N], \ k =2, ..., K_i   \label{eq:f34} \\
                &     &  z_i^{k}        &  \geq   0      &  i & \in [N],\ k \in [K_i] \label{eq:f35}\\ 
                &     & y_{j}           & \in \{ 0, 1\}  &  j & \in [M] \label{eq:f36}
            \end{align} 
            Constraints (\ref{eq:f33}) correspond to Constraints (\ref{ct:zikd}) for $k=1$. Constraints (\ref{eq:f34}) ensure that  $z_i^{k}$ takes the value~1 if $z_i^{k-1}=1$ and if there is no open site at distance $D_i^{k}$ exactly from $i$. Formulations $(F2)$ and $(F3)$ use the same set of variables $y$ and $z$, have exactly the same objective function and describe the same set of optimal integer points with the same linear relaxation~(\cite{Elloumi2010}). However, $(F3)$ has much more zeros in the coefficient matrix, which makes $(F3)$ perform significantly better than $(F2)$.  Therefore, we consider $(F3)$ for our Benders decomposition.

            \subsection{{Resolution methods}\label{subsec:methods}} 
        
                The literature contains many resolution methods for the $(pMP)$. The main heuristics are presented in the following surveys: \cite{Reese2006, Mladenovic2007, Basu2015,  Irawan2015}. For the exact resolution, we only mention the most relevant methods. We refer to~\cite{Marin2019} for more  references. 

                \cite{Galvao1980} solved the $(pMP)$ within a \textit{branch-and-bound} framework solving many linear relaxations of sub-problems of size $N=30$ using formulation $(F1)$. He then devised a method to efficiently obtain good lower bounds instead of optimally solving the relaxed continuous sub-problems.
                
                \cite{avella2007computational} designed a \textit{branch-and-cut-and-price} algorithm also based on $(F1)$ that was able to solve instances up to  $N=5535$. Cuts were added based on new valid inequalities called lifted odd hole and cycle inequalities. Pricing was carried out by solving a master problem to optimality and using dual variables to price out the variables of the initial problem that were not considered in the master, adding new variables if necessary. The novelty of the approach was that Constraints (\ref{eq:f13}) were also relaxed and incorporated to the master problem when the corresponding column was.
                
                \cite{Garcia2011} considered a decomposition algorithm based on formulation $(F2)$. The main idea, also presented  in~\cite{Elloumi2010} on formulation $(F3)$ and implemented in~\cite{elloumi2010computational}, relies on the property that the $z$ variables satisfy $z_i^k \geq z_i^{k+1}$ in any optimal solution of $(F2)$ or its LP relaxation. Therefore, it is enough to solve these problems on a reduced subset of variables $z$, keep enlarging this subset, and stop as soon as one is sure that the remaining $z$ variables can be set to zero to get an optimal solution. This idea is implemented within a \textit{branch-and-cut-and-price} method that the authors name \texttt{Zebra}. It starts  with a very small set of $z$ variables and constraints, and adds more when necessary. \Zebra\;is an exact solution  method that performed well on instances up to $N=85900$ with very large values of $p$. 

                The Benders decomposition has been of great interest in the literature. A survey of this method can be found in~\cite{Rahmaniani2017Literature}. In recent articles, it has produced good results for solving large-scale location problems. For example, \cite{Fischetti2017}  propose an efficient Benders decomposition method within a \BandC\;approach to solve efficiently very large size instances of the \textit{uncapacitated facility location problem} $(UFL)$ where, given a set of potential facility locations and a set of clients, the goal is to find a subset of facility locations to open and to allocate each client to open facilities so that the facility opening plus client allocation costs are minimized. \cite{Cordeau2019benders} described Benders decomposition for two types of location problems:  \textit{the maximal covering location problem} $(MCLP)$, which requires finding a subset of facilities that maximizes the amount of client demand covered while respecting a budget constraint on the cost of the facilities; and \textit{the partial set covering location problem} $(PSCLP)$, which minimizes the cost of the opened facilities while forcing a certain amount of client demand to be covered. They study a decomposition approach of the two problems based on a \textit{branch-and-Benders-cut} reformulation. Their approach is designed for the case in which the number of clients is much larger than the number of potential facility locations.  More recently,  \cite{gaar2021scaleable} perform a Benders decomposition on the \textit{p-Center problem} $(pCP)$. The $(pCP)$  is closely related to the $(pMP)$. The only difference is that instead of minimizing the sum of assignment distances, the largest assignment distance is minimized.

    \section{Benders decomposition for the \textit{(pMP)}\label{sec:Benders} }

        The Benders Decomposition was introduced by~\cite{Benders1962}. The method splits the optimization problem into a \textit{master problem} and one or several \textit{sub-problems}. The master problem and the sub-problems are solved iteratively and at each iteration each sub-problem may add a cut to the master problem.  In this section, we present a Benders decomposition for the \textit{(pMP)} based on  formulation $(F3)$. We show that there is a finite number of Benders cuts and that they can be separated using a polynomial algorithm. 
        
        \subsection{Formulation \label{subsec:descomp}}
     
            For a fixed value of the $y$ variables, the problem decomposes into $N$ sub-problems. Each one computes the assignment distance of a client.
            In the master problem, we remove all $z_i^{k}$ variables and we introduce a new set of continuous  variables $\theta_i$ representing the assignment distance of each client $i\in [N]$:
            \begin{align}
                \nonumber &       &  \displaystyle \min \sum\limits_{i=1}^{N} \theta_i \\
                \nonumber &(MP):  & \sum\limits_{j=1}^{M} y_j   & =  p    &      &  \\
                \nonumber &       & \theta_i                      & \mbox{ satisfies } BD_i & i & \in [N]\\
                \nonumber &       & y_{j}        & \in  \{ 0, 1\}    &    j &\in [M]
            \end{align}

            \noindent where $BD_i$ is the set of benders cuts associated to client $i$. This set is initially empty and grows through the iterations.
            
            The sub-problem for each client  $i \in [N]$ associated to a feasible solution $\bar{y}$ of $(MP)$ is defined by:
            \begin{align}
                \nonumber &       & \displaystyle \min     D_i^1  + \sum\limits_{k=1}^{K_i-1} &(D_i^{k+1}- D_i^{k})z_i^k   &&\\
                \nonumber &(SP_i(\bar{y})):& z_i^1  &  \geq    1 -  \sum\limits_{j:d_{ij} = D_i^1} \bar{y}_{j}  &   \\
                \nonumber &       & z_i^k - z_i^{k-1}    &  \geq  -  \sum\limits_{j:d_{ij} = D_i^k} \bar{y}_{j}  &   k \in\{2, ..., K_i\}  \\
                \nonumber &       &  z_i^{k}             &  \geq   0    &   k \in [K_i] &
            \end{align}
                
            and its corresponding dual sub-problem is:
            \begin{align}
                \nonumber           &&\displaystyle \max     \  D_i^1 +   v^1_{i}(1 -  \sum\limits_{j:d_{ij}=D_i^1}&\bar{y}_j)  -   \sum\limits_{k=2}^{K_i}   v^k_i  \sum\limits_{j:d_{ij}=D_i^k}\bar{y}_j  && \\
                 \nonumber (DSP_i(\bar{y})): && \; v^k_i - v^{k+1}_i & \leq  D_i^{k+1} - D_i^k & k & \in [K_i-1]\\
                \nonumber && \; v^k_i &\geq  0         &    k &\in [K_i] 
            \end{align}

             Note that $(SP_i(\bar{y}))$ and $(DSP_i(\bar{y}))$ are feasible for any $\bar{y}$. From an extreme point $\bar{v}$ of $(DSP_i(\bar{y}))$, we obtain the following optimality Benders cut:
            \begin{equation}
               \theta_i  \geq   D_i^1 +  \bar{v}^1_i(1 -  \sum\limits_{j:d_{ij}=D_i^1}{y}_j)  -  \sum\limits_{k=2}^{K_i}   \bar{v}^k_i \sum\limits_{j:d_{ij}=D_i^k}{y}_j
               \label{cutTHETAi} 
            \end{equation} 

        \subsection{Separation algorithm \label{subsec:separa}}

            The performance of Benders decomposition lies on the resolution of the master problem and the sub-problems. In our decomposition we can have a large number of sub-problems to solve since it is equal to the number of clients at each iteration. Below, we show that the sub-problems can be solved efficiently.
            
            \begin{lemma} \label{prop:solSPi}
               Given a solution $\bar{y}$ of the master problem $(MP)$ or of its LP-relaxation. The following solution  $\bar{z}_i$ is optimal for $(SP_i(\bar{y}))$:
                
                \begin{equation*}
                    \bar{z}_i^{k} =  \max_{k \in [K_i]} \big(0, \; 1- \sum_{j:d_{ij} \leq D_i^{k}} \bar{y}_{j}\big)  \;\;\;  i \in [N]    \label{pr:zik} 
                \end{equation*}
            \end{lemma}
            
            \noindent\textbf{Proof}. We can rewrite the constraints of $(SP_i(\bar{y}))$ for each $i \in [N]$ as follows:
            
            \begin{itemize}
                \item $  z_i^1  \geq 1- \sum\limits_{j:d_{ij} = D_i^{1}} \bar{y}_{j}$ $;$
                \item $  z_i^2 \geq  z_i^1 - \sum\limits_{j:d_{ij} = D_i^{2} }\bar{y}_{j}  \geq 1- \sum\limits_{j:d_{ij} = D_i^{1} } \bar{y}_{j}  -\sum\limits_{j:d_{ij} = D_i^{2}}  \bar{y}_{j}  = 1 -\sum\limits_{j:d_{ij} \leq D_i^{2}}  \bar{y}_{j} $ $;$
                \item ...
                \item $  z_i^k \geq  z_i^{k-1} - \sum\limits_{j:d_{ij} = D_i^{k} }\bar{y}_{j}  \geq 1- \sum\limits_{j:d_{ij} = D_i^{k-1} } \bar{y}_{j}  -\sum\limits_{j:d_{ij} = D_i^{k}}  \bar{y}_{j}  = 1 -\sum\limits_{j:d_{ij} \leq D_i^{k}}  \bar{y}_{j}$ $.$
            \end{itemize}
            \vspace{0.2cm}

            Since we minimize an objective function with non-negative coefficients, the $\bar{z}_i^k$ variables are as small as possible in an optimal solution. Thus, setting $\bar{z}_i^k$ to the value $\max\big(0, \; 1- \sum_{j:d_{ij} \leq D_i^{k}} \bar{y}_{j}\big)$ leads to a feasible and optimal solution.
            
            $ \hfill\square$

            From Lemma~\ref{prop:solSPi} we observe that the optimal values of variables $z_i^k$ in $(SP_i(\bar y))$ are decreasing when $k$ increases. In order to obtain a dual solution, we identify the last strictly positive term of this sequence.

            \begin{definition}
              Given a solution $\bar{y}$ of the master problem $(MP)$ or of its LP-relaxation. Let $ \Tilde{k}_i $ be the following index,
                
                \begin{equation*}
                    \Tilde{k}_i =
                    \begin{dcases}
                        0   & \text{if}  \quad   \sum_{j:d_{ij} = D_i^{1}} \bar{y}_{j} \geq 1 \\
                         \max\{k\in [K_i] :  1- \sum_{j:d_{ij} \leq D_i^{k}} \bar{y}_{j} >0 \} & \text{otherwise}   
                    \end{dcases}
                      \hspace{1cm} i \in [N]  \label{ct:ktildei}
                \end{equation*}
            \end{definition}

            Note that, if $\bar{y}$ is binary, then the assignment distance of client $i$  in the feasible solution $\bar y$ is $D_i^{\Tilde k_i+1}$. 
        
            \begin{lemma} \label{prop:optSPi}
                Given a solution $\bar{y}$ of the master problem $(MP)$ or of its LP-relaxation and the corresponding indices $\tilde{k}_i$,  the optimal value of $SP_i(\bar{y})$ for $i \in [N]$ is:
                
                \begin{equation}
                    OPT(SP_i(\bar{y}))  =
                    \begin{dcases}
                     D_i^{1} & \text{if     } \tilde{k}_i =0 \\
                     \; \\
                    D_i^{\Tilde{k}_i+1}  -    
                 \sum\limits_{j:d_{ij}\leq D_i^{\tilde{k}_i}}(D_i^{\Tilde{k}_i+1}- d_{ij}) \bar{y}_j   & \text{otherwise}   
                    \end{dcases}
                \end{equation} 
            \end{lemma}        
                
            \noindent\textbf{Proof}.  If  $\tilde{k}_i =0$, then using Lemma~\ref{prop:solSPi}  and Definition~\ref{ct:ktildei},then all the values $\bar z_i^k$ are equal to 0 and $OPT(SP_i(\bar y))=0$. Otherwise, when $ \tilde{k}_i \geq 1$, we have the following:

            \begin{align*}
                OPT(SP_i(\bar{y}))  &=   D_i^1  + \sum\limits_{k=1}^{\tilde{k}_i} (D_i^{k+1}- D_i^{k})   (1- \sum_{j:d_{ij} \leq D_i^{k}} \bar{y}_{j})  \\
                & =  D_i^1 +  D_i^{\tilde{k}_i+1} (1- \sum_{j:d_{ij} \leq D_i^{\tilde k_i}} \bar{y}_{j}) -D_i^1 (1- \sum_{j:d_{ij} \leq D_i^1} \bar{y}_{j} ) + \\ & \quad \quad  \sum\limits_{k=2}^{\tilde{k}_i} D_i^{k}   ((1- \sum_{j:d_{ij} \leq D_i^{k-1}} \bar{y}_{j})  -(1- \sum_{j:d_{ij} \leq D_i^{k}} \bar{y}_{j} )) \\
                & =  D_i^{\tilde{k}_i+1} (1- \sum_{j:d_{ij} \leq D_i^{\tilde k_i}} \bar{y}_{j}) + \sum\limits_{k=1}^{\tilde{k}_i}  \sum_{j:d_{ij} = D_i^{k}} d_{ij} \bar{y}_{j}  \\
                & =  D_i^{\tilde{k}_i+1} - \sum_{j:d_{ij} \leq D_i^{\tilde k_i}} (D_i^{\tilde{k}_i+1}-  d_{ij}) \bar{y}_{j}
            \end{align*} \hfill $\square$
                
            \begin{proposition}\label{prop:optDSPi}
            Given a solution $\bar{y}$ of the master problem $(MP)$ or of its LP-relaxation and the corresponding indices $\tilde{k}_i$. The following solution $\bar{v}_i$ is optimal for $DSP_i(\bar{y})$:
                    \begin{equation*}
                            \bar{v}^k_i = 
                        \begin{dcases}
                            D_i^{\Tilde{k}_i+1}-D_i^{k},              & \mbox{if } k \leq \Tilde{k}_i \\
                            0,                    & \mbox{otherwise}\\
                        \end{dcases}
                        \hspace{1cm} \ i \in [N] \ k \in [K_i]  
                     \end{equation*} 
            \end{proposition}

            \noindent\textbf{Proof}. From the theorem of complementary slackness,  for $i \in [N]$ we have:
           
            \begin{equation}
                \bar{v}^k_i-\bar{v}^{k+1}_i = D_i^{k+1} - D_i^{k} \hspace{1cm}   k \leq \Tilde{k}_i \label{eqviks}
            \end{equation}
            because for $k\leq \tilde k_i$, the primal variables values $\bar z_i^k$ are strictly positive. Then, for $k \in [\tilde{k}_i]$ if we sum  Equations (\ref{eqviks}) from $k'=k$ to $k'=\tilde{k}_i$,  we obtain 
            
            \begin{equation}
                    \bar{v}^{k}_i = v_{i}^{\tilde{k}_i+1} + (D_i^{\tilde{k}_i+1} - D_i^{k})  \hspace{1cm}   k \leq \Tilde{k}_i \label{eqviks2}
            \end{equation}
            We build a feasible solution of $(DSP_i(\bar{y}))$ by setting $\bar{v}^k_i=0$ for  $k > \tilde{k}_i$ and   $\bar{v}^k_i = (D_i^{\tilde{k}_i+1} - D_i^{k}) $ for $k \leq \Tilde{k}_i$. The objective value $\mathcal V(\bar v)$  of this solution  is:
        
            \begin{align*}
                \mathcal V(\bar v)  &= D_i^1 +   (D_i^{\tilde{k}_i+1} -     D_i^{1})(1 -  \sum\limits_{j:d_{ij}=D_i^1}\bar{y}_j)  -   \sum\limits_{k=2}^{\tilde{k}_i}   (D_i^{\tilde{k}_i+1} - D_i^{k})  \sum\limits_{j:d_{ij}=D_i^k}\bar{y}_j  \\
                &=D_i^{\tilde{k}_i+1}  -   \sum\limits_{k=1}^{\tilde{k}_i}   (D_i^{\tilde{k}_i+1} - D_i^{k})  \sum\limits_{j:d_{ij}=D_i^k}\bar{y}_j  \\
                &=   D_i^{\tilde{k}_i+1}  - \sum\limits_{k=1}^{\tilde{k}_i} \sum_{j:d_{ij} = D_i^{k}}(D_i^{k+1}- D_i^{k}) \bar{y}_{j} \\ 
                & =  D_i^{\Tilde{k}_i+1}  -    
                \sum\limits_{j:d_{ij}\leq D_i^{\tilde{k}_i}}(D_i^{\Tilde{k}_i+1}- d_{ij}) \bar{y}_j   = OPT(SP_i(\bar y))
                \label{opt:dspi}       
             \end{align*}
                 
            Hence, this dual solution is feasible and it has the same objective value as the optimal value of the primal solution. Thus, this solution $\bar v$ is optimal for $(DSP_i(\bar y))$.
    
            \hfill $\square$
            
            \noindent\textbf{Corollary~1} The Benders cuts (\ref{cutTHETAi}) can be written as follows:
            \begin{equation}
                \begin{dcases}    \theta_i  \geq D_i^{1} &  \mbox{ if }\Tilde{k}_i =0 \\
                    \; \\
                    \theta_i \geq D_i^{\Tilde{k}_i+1}  -    
                 \sum\limits_{j:d_{ij}\leq D_i^{\tilde{k}_i}}(D_i^{\Tilde{k}_i+1}- d_{ij}) {y}_j &  \mbox{otherwise}
               \label{cutTHETAi5} 
               \end{dcases}
            \end{equation} 
                
            This corollary highlights that the Benders cuts can be obtained in polynomial time by computing $\tilde{k}_i$. We use Algorithm~\ref{alg:Separation} to separate Constraints~\eqref{cutTHETAi5}. For each client $i\in[N]$, we  first compute $\Tilde k_i$ and $OPT(DSP_i(\bar{y}))$ from the current $(MP)$ solution $(\bar{y}, \ {\bar{\theta}})$ (steps~3 and~4). Then, if the value of the assignment distance in the current $(MP)$ solution is underestimated (step~6) we directly construct the corresponding Benders cuts (\ref{cutTHETAi5}) (step~7).

            Let  $\of_{ir}\in [M]$ be the $r^{th}$ closest site to client $i$. Hence, $d_{i\of_{ir}}$ is the distance from $i$ to its $r^{th}$  closest site. 
            The  complexity of Algorithm~\ref{alg:Separation} is determined by the computation of index $\tilde{k}_i$. We compute it in $O(M)$ by Algorithm~\ref{alg:kiindex} which uses $\of_{ir}$ to obtain indices $\tilde k_i$ satisfying Definition 1. Observe that, the $N\times M$ matrix $\of$ can be built only once in a preprocessing step. 
            
            Then, given $\tilde{k}_i$ it is easy  to compute step 4 and 5 of Algorithm~\ref{alg:Separation} in $O(M)$ and $O(1)$ respectively. Consequently, considering the $N$ clients, a complexity in $O(NM)$ is obtained for Algorithm~\ref{alg:Separation}.

            \resizebox{0.9\textwidth}{!}{%
            \centering
            \begin{algorithm} [H]\label{alg:Separation}
                \caption{Separation algorithm}
                \DontPrintSemicolon
                \SetKwInOut{Input}{input}
                \SetKwInOut{Output}{output}
                \SetKwInOut{Init}{init}
                \Input{ 
                    \begin{itemize}[noitemsep,topsep=0pt]
                        \item  Instance data $([N]$, $[M]$, $[K_i]$, distances ${D_i^0, ... , D_i^{K_i}}$ and $d_{ij}$ for each $i \in [N]$, $j \in [M])$
                        \item  Current $(MP)$ solution $(\bar{y}$,  \ $\bar\theta)$
                    \end{itemize}
                }    
                \Output{                
                    \begin{itemize}[noitemsep,topsep=0pt]
                         \item Upper bound of (MP).
                    \end{itemize}
                }
                \vspace{0.3cm}
                $UB\leftarrow 0$\;
                \For{$i \in [N]$}{
                     Compute $\tilde{k}_i$ with Algorithm~\ref{alg:kiindex}  \;
                     Compute $OPT(DSP_i(\bar{y}))$\ (more precisely $OPT(SP_i(\bar y))$) as in Lemma~\ref{prop:optSPi} \;
                     $UB \leftarrow UB + OPT(DSP_i(\bar{y}))$\;
                     \If{ $\bar\theta_i < OPT(DSP_i(\bar{y}))$}{  
                       Add the corresponding cut (\ref{cutTHETAi5}) to (MP)
                     }
                }
                \Return UB\;
            \end{algorithm}
            }

            \vspace{0.5cm}
            \resizebox{0.9\textwidth}{!}{%
            \centering
            \begin{algorithm} [H]\label{alg:kiindex}
                \caption{Computing $\Tilde{k}_i$ }
                \SetKwInOut{Input}{input}
                \SetKwInOut{Output}{output}
                \SetKwInOut{Init}{init}
                \DontPrintSemicolon
                \Input{ 
                    \begin{itemize}[noitemsep,topsep=0pt]
                        \item  Instance data $([N]$, $[M]$, $[K_i]$, $\of$ matrix, and distances $d_{ij}$  for $i \in [N]$, $j \in [M])$
                        \item  Current $(MP)$ solution $\bar{y}$
                        \item $i\in[N]$
                    \end{itemize}
                }    
                \vspace{0.2cm}
                
                \Output{ \begin{itemize}[noitemsep,topsep=0pt]
                            \item  The index $\Tilde{k}_i$   associated to $\bar{y}$ as defined in Definition~\ref{ct:ktildei}
                         \end{itemize}
                }
                \vspace{0.3cm}
                $\Tilde{k}_i \leftarrow 0 $ \;
                $r  \leftarrow 1$ \;
                $ val \leftarrow    1 - \bar{y}_{\of_{ir}}$ \;
                
                \While{$val >0$ and $r < M$}{
                \If{$d_{i(\of_{i(r+1)})} > d_{i\of_{ir}} $}{
                $\Tilde{k}_i  \leftarrow \Tilde{k}_i + 1 $ \; 
                }
                $ r \leftarrow (r +1)$ \;
                $val \leftarrow    val - \bar{y}_{\of_{ir}}$ \;     
                }
                \Return $\Tilde k_i$\;
                \end{algorithm}
            }
            
        \subsection{A resulting new formulation}\label{subsec:newform}
         
            The Bender cuts~\eqref{cutTHETAi5} lead to the following new compact formulation for $(pMP)$: 
            \begin{align} 
                \nonumber &       &  \displaystyle \min \sum\limits_{i=1}^{N} \theta_i \\
                \nonumber &(F4):  & \sum\limits_{j=1}^{M} y_j   & =  p    &      &  \\
                 &       & \theta_i                      & \geq     D_i^{{k}+1}  - \sum_{j:d_{ij} \leq D_i^{{k}}}(D_i^{{k}+1}- d_{ij}) y_{j}  &    i &\in [N],\ k \in [K_i]\; \label{eq:bdcuts}\\ 
                \nonumber &       & y_{j}        & \in  \{ 0, 1\}    &    j &\in [M]
            \end{align}
            
            Constraints~(\ref{eq:bdcuts}) ensure that each variable $\theta_i$ is larger than $D_i^{{k}+1}$ unless a site is open at a smaller distance than $D_i^{{k}}$ from $i$. This formulation $(F4)$  has $(N+M)$ variables which is less  than  $(F2)$ and $(F3)$ but it has the same number of constraints. Nevertheless, the constraint matrix is roughly as dense as $(F2)$ and it has the same continuous relaxation. 
            
            Table~\ref{tab:comparison} presents the results of four formulations of the $(pMP)$ on five instances described in Section~\ref{subsec:definitions}. A time limit of 600 seconds is considered. For each formulation the relative optimality gap ($gap(\%)$) and the resolution time ($t(s)$) in seconds are presented.

            \begin{table}[H]
                \centering
                \resizebox{0.85\textwidth}{!}{%
                \begin{tabular}{@{}|cccc|cc|cc|cc|cc|@{}}
                \toprule
                \multicolumn{4}{|c}{INSTANCE} & \multicolumn{2}{|c}{F1}                         & \multicolumn{2}{|c}{F2}                         & \multicolumn{2}{|c}{F3}                         & \multicolumn{2}{|c|}{F4}                         \\ \midrule
                Name    & N=M  & p  & OPT             & \multicolumn{1}{c}{$gap(\%)$} & \multicolumn{1}{c|}{$t(s)$} & \multicolumn{1}{c}{$gap(\%)$} & \multicolumn{1}{c|}{$t(s)$} & \multicolumn{1}{c}{$gap(\%)$} & \multicolumn{1}{c|}{$t(s)$} & \multicolumn{1}{c}{$gap(\%)$} & \multicolumn{1}{c|}{$t(s)$} \\ \midrule
                pmed26  & 600  & 5  & \textbf{9917}   & 0\%            & 228          & 0\%            & 40           & 0\%            & 8            & 0\%            & 57           \\
                pmed31  & 700  & 5  & \textbf{10086}  & 0\%            & 282          & 0\%            & 36           & 0\%            & 7            & 0\%            & 58           \\
                pmed35  & 800  & 5  & \textbf{10400}  & 0\%            & 527          & 0\%            & 104          & 0\%            & 9            & 0\%            & 95           \\
                pmed38  & 900  & 5  & \textbf{11060}  & 74,1\%           & 600          & 0\%            & 75           & 0\%            & 19           & 0\%            & 115          \\
                pmed39  & 900  & 5  & \textbf{11069}  & 10,7\%           & 600          & 0\%            & 66           & 0\%            & 19           & 0\%            & 105          \\
                pmed40  & 900  & 5  & \textbf{12305}  & 0\%            & 579          & 0\%            & 60           & 0\%  & 10           & 0\%            & 104          \\ \bottomrule
                \end{tabular}
                }
                \caption{Comparison between different (pMP) formulations. Time limit of 600 seconds.}
                \label{tab:comparison}
            \end{table}

            Results in Table~\ref{tab:comparison} confirm the expected performance between formulations $(F1)$, $(F2)$ and $(F3)$ already described in Section~\ref{subsec:formulations}. Moreover, we see that $(F4)$ takes more time than $(F2)$ and $(F3)$. 

        \subsection{Decomposition algorithm implementation \label{subsec:algoBenders}}
        
            To improve the performances, we implement a two-phase Benders decomposition algorithm. Let $(\overline{MP})$ be the master problem without the integrity constraints. One classic implementation is to first solve Benders decomposition for $(\overline{MP})$ (\emph{Phase~1}). Then, both the Benders cuts generated in the first phase and the integrity constraints are added and the obtained master problem $(MP)$ is solved through a \textit{branch-and-Benders-cut} (\emph{Phase~2}).

            \subsubsection{Phase~1: Solving the linear relaxation of the  master problem}
            Phase~1 is summarized in Algorithm~\ref{alg:SolveMp}. The current master problem $\overline{MP}$ is solved at step 4 through a linear programming solver and provides a candidate solution $(\bar y, \overline\theta)$ while the sub-problems are solved at steps 2 and 6 using Algorithm~\ref{alg:Separation}.  To enhance the performance this phase includes the following improvements:
                
            \begin{itemize}
                \item \textbf{Initial solution}: Providing a good candidate solution to the initial $(\overline{MP})$ can significantly reduce the number of iterations. Consequently, as~\cite{Garcia2011} we compute a first solution using the \Popstar\;heuristic~(\cite{Resende2004}) which, to the best of our knowledge, is the best current heuristic for the $(pMP)$. \Popstar\;is a hybrid heuristic that combines elements of several metaheuristics. It uses a multi-start method in which a solution is built at each iteration as in a GRASP algorithm. It is followed by an intensification strategy, with a tabu search and a scatter search. And in a post-optimization phase, they use the concept of multiple generations, a characteristic of genetic algorithms. The solution $y^h$ provided by this heuristic and its objective value $UB^h$ are inputs of Algorithm~\ref{alg:SolveMp}. The latter is used to initialize $UB^1$ at step $1$.
                \item \textbf{Rounding heuristic}: Since in Phase~1 most of the solutions provided by $(\overline{MP})$ are fractional, we use a primal heuristic to try to improve the upper bound of the problem. At each iteration we open the sites associated to the $p$ largest values of $\bar{y}$ (steps 7 to 11 in Algorithm~\ref{alg:SolveMp}).
            \end{itemize}

            The objective value of $(\overline{MP})$ and the  sub-problem optimal value $OPT(SP)$ allow us to update the lower bound $LB_{\overline{MP}}$ (step 5) and the upper bound $UB_{\overline{MP}}$ (step 8), respectively. In each iteration the rounding heuristic tries the improve $UB^1$ (step 10). The iterative algorithm is terminated by updating $LB^1$ when $LB_{\overline{MP}} = UB_{\overline{MP}}$.
                   
            \subsubsection{Phase~2:  Solving the master problem with \textit{branch-and-Benders-cut} approach}

                Once the  continuous relaxation of $(MP)$ is solved by Phase~1, we keep the generated cuts and add the integrity constraints on variables $y$. We use a \BandC~algorithm in which we solve the sub-problems at each node which provides an integer solution in order to generate Benders cuts. This approach is called \textit{branch-and-Benders-cut}~(\cite{Rahmaniani2017Literature}). The resolution of the sub-problems is performed through callbacks which is a feature provided by mixed integer programming solvers. In order to enhance the performance of Phase~2,  we implement the following improvements:
                  
             \begin{itemize}
                \item \textbf{Constraint reduction}:
                At the end of Phase~1, most of the generated Benders cuts  are not saturated by the current fractional solution. We remove most of them to reduce the problem size. The cuts of a client $i$ are related to indexes $\tilde{k}_i$ obtained at different iterations. Let $\hat{k}$ be the highest of these indexes associated with a saturated constraint. We remove all constraints of client $i$ which associated index is higher than $\hat{k}$. This reduction performs better than removing all unsaturated constraints.
                \item \textbf{Reduced cost fixing}:
                At the end of Phase~1, given the bounds $LB^1$ and $UB^1$, we can perform an analysis of the reduced costs $\overline{rc}$ of the last fractional solution $\overline y$ provided by Algorithm~\ref{alg:SolveMp}. For any site $j$ such that  $LB^1 + \overline{rc}_j > UB^1$, $y_j$ can be set to 0. Similarly, for any site $j$ such that $LB^1 - \overline{rc}_j > UB^1$, $y_j$ can be fixed to 1. We computationally observed that this rule is efficient  in instances where $p$ is small.
             \end{itemize}

                \resizebox{0.9\textwidth}{!}{%
                \centering
                \begin{algorithm}[H]\label{alg:SolveMp}
                \caption{Phase~1 - Solving $(\overline{MP})$}
                \DontPrintSemicolon
                \SetKwInOut{Input}{input}
                \SetKwInOut{Output}{output}
                \SetKwInOut{Init}{init}
                \Input{\justifying  
                    \begin{itemize}[noitemsep,topsep=0pt]
                        \item Instance data ($N$, $M$, $p$ , Distances $D_i^0, ... , D_i^{K_i}$ and $d_{ij}$ with $i \in [N]$, $j \in [M]$)
                        \item Heuristic $(pMP)$ solution $y^h$ with value $UB^h$
                    \end{itemize} 
                }
                \Output{\justifying
                    \begin{itemize}[noitemsep,topsep=0pt]
                        \item Lower bound $LB^1$ and a feasible integer solution ${y^1}$ with value $UB^1$
                    \end{itemize} 
                }
                \vspace{0.2cm}
                     
                $(LB_{\overline{MP}}, UB_{\overline{MP}}, y^1, UB^1) \leftarrow (0, UB^h, y^h, UB^h)$ \\
                Use Algorithm 1 to add the Benders cuts associated with $y^h$  to $(\overline{MP})$ \\
                  
                \vspace{0.2cm}
                \While{$LB_{\overline{MP}} < UB_{\overline{MP}}$} {
                     \vspace{0.2cm}
                    $(\bar{y}$, ${\overline\theta})\leftarrow $ Solution of $(\overline{MP})$\;
                    $ LB_{\overline{MP}} \leftarrow  \sum\limits_{i=1}^N {\overline\theta}_i $\\  
                    Use Algorithm 1 to compute the sub-problem optimal value $OPT(SP)$ and add Benders cuts associated with $\bar{y}$ to $(\overline{MP})$\;
                    \If{   $  OPT(SP) < UB_{\overline{MP}} $}{
                       $ UB_{\overline{MP}} \leftarrow  OPT(SP)$ 
                     }
                    \If{$\bar{y}$ is fractional}{
                        $({y}^r, UB^r)\leftarrow$  Rounded solution from $\bar{y}$ and its value\;
                      \If{$UB^r < UB^1$}{  
                         $UB^1 \leftarrow UB^r$\;
                         ${y^1} \leftarrow {y}^r$\;
                      }
                    }
                }
                $LB^1\leftarrow LB_{\overline{MP}}$\;
                \Return $LB^1, {y^1},UB^1$\;
            \end{algorithm} 
                }

    \section{Computational study \label{sec:exper}}    
    
        In this section, we compare the results of our Benders decomposition method with those of the  state-of-the-art methods described in Section~\ref{subsec:methods}.
        
        \subsection{Instances and technical specifications \label{subsec:definitions}}
        
        We consider the same instances  used in~\cite{Garcia2011} that is the  $p$-median instances from OR-Library~(\cite{ORlib}) and TSP instances from the TSP-Library~(\cite{TSPlib}). In all these instances, the sites are at the same location as the clients and thus $N=M$. The set of OR-Library contains instances with 100 to 900 clients, and the value of $p$ is between 5 and 500. The set of TSP-Library selected contains between 1304 and 115475 clients. All customer demand points are given as two-dimensional coordinates, and the Euclidean distance rounded down to the nearest integer is used as distance. This follows previous literature.

        We also consider the RW instances originally proposed by~\cite{Resende2004} with the \Popstar\;heuristic. They correspond to completely random distance matrices. The distance between each site and each client is an integer value taken uniformly in the interval $[1,n]$. Moreover, the distance between client $i$ and site $j$ is not necessarily equal to the distance between site $j$ and client $i$. Four different values of $N=M$ are considered: 100, 250, 500, and 1000.
        
        Finally, we include in our experimentation the ODM instances which were introduced by~\cite{briant2004optimal} and are used by~\cite{avella2007computational} with a \textit{branch-and-cut-and-price} algorithm. We name this algorithm \Avella. These instances correspond to the \textit{Optimal Diversity Management Problem} which can be treated as a $p$-median problem on a graph. For this problem there exist instances with $N$ equals to 1284, 33773, and 5335. It was already observed by~\cite{avella2007computational} that instances with $N$ equal to 1284 and 5335 are easy to solve.  Therefore, we have only considered the instances with the value of $N=3773$.
        
        Our study was carried out on an Intel XEON W-2145 processor 3,7 GHz, with 16 threads, but only 1 was used, and 256 GB RAM. IBM ILOG CPLEX 20.1 was used as \textit{branch-and-cut} framework. We apply the described separation algorithm in the GenericCallback of CPLEX, which gets called whenever a feasible integer solution is found. We set optimality tolerance EpGap to $10^{-10}$, the absolute tolerance on the gap between the best integer objective and EpAGap to 0.9999, the objective of the best remaining node, also called absolute MIP gap. Considering that our Benders decomposition can easily find feasible solutions, we have set the MIP emphasis switch to BestBound in order to prove the optimality as fast as possible. We  set the branch-up-first parameter BRDIR  to 1, since this tends to produce branching trees with fewer nodes. We use a time limit of 36000 seconds for Phase~2, indicated by TL in the tables when this is reached.
        
        We were able to run the \Zebra\;and \Popstar\;methods on our computer. \Zebra\;code was provided by the authors of~\cite{Garcia2011}  and \Popstar\;code is available online\footnote{ http://mauricio.resende.info/popstar/}. However, we cannot report an up-to-date time for \Avella\;algorithm which was originally carried out on a Compaq EVO W4000 Personal Computer with Pentium IV-1.8 Ghz processor and 1 Gb RAM using the LP solver IBM ILOG CPLEX 8.0 with a time limit of 100 hours per instance (indicated by TL2 in the corresponding table).

        \subsection{Performance analysis \label{subsec:perform}}

            The results for the different instances are presented below. The information in the tables is organized as follows:
      
            \begin{itemize}
            \item Instance data
                \begin{itemize}
                    \item name: name of the instance.
                    \item $N=M$: size of the instance (number of clients equal to the number of sites).
                    \item $p$: number of sites to open.
                    \item $OPT / BKN$: optimal value of the instance (in \textbf{bold}) if it is known or the best-known solution value obtained given the time limit, otherwise. If the value is \underline{underlined},  it means that it is the first time the instance is solved to optimality or that we improve the best-known value.
                \end{itemize}
            \item Phase~1 results:
                \begin{itemize}
                        \item $LB^1$: lower bound of $(pMP)$ obtained at the end of Phase~1.
                        \item $UB^1$: upper bound of $(pMP)$ obtained at the end of Phase~1.
                        \item $T^1$: CPU time in seconds required to complete Phase~1.
                    \end{itemize}
            \item Phase~1 + Phase~2 results:
                \begin{itemize}
                        \item $gap$: relative optimality gap between the lower and upper bound obtained at the end of Phase~2.
                        \item $iter$: number of total iterations required for the Benders decomposition, i.e., the number of times a fractional solution or an integer solution was separated in Phase~1 and Phase~2, respectively.
                        \item $T^{tot}$: the total CPU time in seconds required to exactly solve the instance. 
                \end{itemize}
            \item \Zebra,  \Popstar\;and \Avella\;results:
                    \begin{itemize}
                        \item $gap$:  relative optimality gap between the lower and upper bound obtained at the end of the  method.
                        \item $T$: total CPU time in seconds required to complete the algorithms of \cite{Garcia2011}),    \cite{Resende2004} or \cite{avella2007computational}. A diamond ($\blacklozenge$) means that the computer ran out of memory while solving the problem. 
                \end{itemize}
                \item Average total time 
                    \begin{itemize}
                        \item the average total time by our method and \Zebra\;is presented in the corresponding tables. This average is calculated considering only the instances where both methods solve the instances to optimality.
                    \end{itemize}
            \end{itemize}

            \subsubsection{ORLIB Instances}
        
                ORLIB instances are considered to be easy instances. \cite{Garcia2011} only include in their results the most difficult ones. The results are presented in Table~\ref{tab:ORLIB}. In all instances we reached the optimal value within a few seconds. The computation time with respect to~\cite{Garcia2011} is significantly reduced in the instances with a small value of $p$. When $p$ is large, \Zebra\;can be faster for a few instances solved in less than 1 second.
        
            \subsubsection{TSP Instances}
            
                The results on TSP instances are presented in Tables ~\ref{tab:TSPmedium},~\ref{tab:TSPlarge}, and~\ref{tab:TSPhuge}, for medium, large and huge instances, respectively. Our method reaches the optimal solution in most instances. Very good  $LB^1$ and $UB^1$ bounds are quickly found at the end of Phase~1. 
                
                {In Tables~\ref{tab:TSPmedium} and~\ref{tab:TSPlarge}} we can observe that 10 medium and large instances are not solved optimally by \Zebra\;due to a lack of memory or for reaching the time limit. 
                However, our method does not face any memory problem and only 2 medium and 2 large instances reach the time limit of 10 hours with an optimality gap lower than $0,1\%$. 
                   
                Regarding the huge instances in Table~\ref{tab:TSPhuge},  we solve 38 out of the 44 instances whereas \Zebra\;only solves 16 instances due to a lack of memory.  For all the huge instances, the rounding heuristic step of Phase~1, the reduced cost fixing step, and the constraint reduction step of Phase~2 were not used as they take too much time. Nevertheless, we use the rounding heuristic once at the end of Phase~1 to update $UB^1$. For the largest instances ch71009, pla85900, and usa15475, we use  a randomly generated solution instead of \Popstar\;which takes a long time to compute a solution for these instances. With our method, the time limit is reached in only 7 instances, in which the optimality gap is lower than $0,9\%$.  We can also provide for the first time the optimal values for 7  instances with $N=115455$.
            
            \newpage
            \subsubsection{RW instances}
    
                The results on RW instances are summarized in Table~\ref{tab:RW}. Even small RW instances can be very difficult to solve as previously observed by~\cite{elloumi2010computational}. We think that it is mainly due to the fact that the instances are non-euclidean. Furthermore, the total number of distances $K$ is closer  to $N\times M$, leading to more variables and constraints in Formulation $(F_2)$. For large values of $p$, our decomposition can quickly solve the instances to optimality.  These instances were not considered by either~\cite{avella2007computational} or~\cite{Garcia2011}. Moreover, the code of \Zebra\;cannot handle non-symmetric instances.  Consequently,  we only report the computation time and value $UB^h$ of the solution computed by the heuristic \Popstar. 
        
            \subsubsection{ODM instances}
        
                The results on ODM instances are summarized in Table~\ref{tab:ODM}. To solve these instances we need to add $N$ constraints to ensure that each client is assigned to one of its non-forbidden neighbors. This generates a more complex master problem to solve. To save computation time, we did not use the rounding heuristic of Phase~1.
                
                Obtaining a solution for these instances is hard since \Popstar\;does not support their format and since random solutions may not be feasible due to the sparsity of the graphs. Consequently, to obtain an initial solution, we initiate a resolution of its corresponding formulation $(F3)$ by CPLEX and stop it once it has found 3 feasible solutions. This approach empirically proved to be a good compromise between computation time and optimality gap. 
                We were also unable to use \Zebra\; for these instances. We solve to optimality all these instances.

                \begin{table}[H]
                    \centering                      \resizebox{0.62\textwidth}{!}
                    {%
                    \begin{tabular}{@{}lccr|rrr|rcr|rr@{}}
                    \toprule
                    \multicolumn{4}{c|}{INSTANCE}                                                                   & \multicolumn{3}{c|}{PHASE 1}                                         & \multicolumn{3}{c|}{PHASE 1 + 2}                                                                                      & \multicolumn{2}{c}{Zebra}                                                             \\ \midrule
                    \multicolumn{1}{c}{Name} & \multicolumn{1}{c}{$N=M$} & \multicolumn{1}{c}{$p$} & \multicolumn{1}{c|}{\begin{tabular}[c]{@{}c@{}}$OPT /$\\ $BKN$\end{tabular}} & \multicolumn{1}{c}{$LB^1$} & \multicolumn{1}{c}{$UB^1$} & \multicolumn{1}{c|}{$T^1$} & \multicolumn{1}{c}{$gap$} & \multicolumn{1}{c}{$iter$} & \multicolumn{1}{c|}{$T^{tot}$} & \multicolumn{1}{c}{$gap$} & \multicolumn{1}{c}{T} \\ \midrule
                    pmed26                     & 600                      & 5                      & \textbf{9917}            & 9854                     & 9917                     & 0,14                    & 0\%                                                                    & 15                        & 1,12                      & 0\%                                                                    & 5,58                   \\
                    pmed27                     & 600                      & 10                     & \textbf{8306}            & 8302                     & 8307                     & 0,18                    & 0\%                                                                    & 19                        & 0,63                      & 0\%                                                                    & 0,91                   \\
                    pmed28                     & 600                      & 60                     & \textbf{4498}            & 4498                     & 4498                     & 0,12                    & 0\%                                                                    & 7                         & 0,16                      & 0\%                                                                    & 0,12                   \\
                    pmed29                     & 600                      & 120                    & \textbf{3033}            & 3033                     & 3033                     & 0,08                    & 0\%                                                                    & 8                         & 0,11                      & 0\%                                                                    & 0,05                   \\
                    pmed30                     & 600                      & 200                    & \textbf{1989}            & 1989                     & 1989                     & 0,04                    & 0\%                                                                    & 9                         & 0,04                      & 0\%                                                                    & 0,04                   \\
                    pmed31                     & 700                      & 5                      & \textbf{10086}           & 10026                    & 10086                    & 0,16                    & 0\%                                                                    & 13                        & 0,93                      & 0\%                                                                    & 4,23                   \\
                    pmed32                     & 700                      & 10                     & \textbf{9297}            & 9293                     & 9297                     & 0,27                    & 0\%                                                                    & 9                         & 0,91                       & 0\%                                                                    & 1,29                   \\
                    pmed33                     & 700                      & 70                     & \textbf{4700}            & 4700                     & 4700                     & 0,18                    & 0\%                                                                    & 8                         & 0,24                      & 0\%                                                                    & 0,17                   \\
                    pmed34                     & 700                      & 140                    & \textbf{3013}            & 3013                     & 3013                     & 0,08                    & 0\%                                                                    & 7                         & 0,08                      & 0\%                                                                    & 0,07                   \\
                    pmed35                     & 800                      & 5                      & \textbf{10400}           & 10302                    & 10400                    & 0,17                    & 0\%                                                                    & 9                         & 1,24                      & 0\%                                                                    & 5,32                   \\
                    pmed36                     & 800                      & 10                     & \textbf{9934}            & 9834                     & 9934                     & 0,25                    & 0\%                                                                    & 18                        & 27,2                      & 0\%                                                                    & 26,7                   \\
                    pmed37                     & 800                      & 80                     & \textbf{5057}            & 5057                     & 5057                     & 0,22                    & 0\%                                                                    & 6                         & 0,29                      & 0\%                                                                    & 0,15                   \\ \midrule
                    pmed38                     & 900                      & 5                      & \textbf{11060}           & 10948                    & 11060                    & 0,23                    & 0\%                                                                    & 15                        & 4,42                      & 0\%                                                                    & 10,1                  \\
                    pmed38                     & 900                      & 10                     & \textbf{9431}            & 9362                     & 9431                     & 0,31                    & 0\%                                                                    & 20                        & 4,17                      & 0\%                                                                    & 10,3                   \\
                    pmed38                     & 900                      & 20                     & \textbf{7839}            & 7832                     & 7839                     & 0,32                    & 0\%                                                                    & 11                        & 1,03                      & 0\%                                                                    & 2,87                   \\
                    pmed38                     & 900                      & 50                     & \textbf{5892}            & 5889                     & 5892                     & 0,26                    & 0\%                                                                    & 11                        & 0,86                      & 0\%                                                                    & 1,19                   \\
                    pmed38                     & 900                      & 100                    & \textbf{4450}            & 4450                     & 4450                     & 0,22                    & 0\%                                                                    & 8                         & 0,31                       & 0\%                                                                    & 0,15                   \\
                    pmed38                     & 900                      & 200                    & \textbf{2905}            & 2905                     & 2905                     & 0,09                    & 0\%                                                                    & 8                         & 0,09                      & 0\%                                                                    & 0,08                   \\
                    pmed38                     & 900                      & 300                    & \textbf{1972}            & 1972                     & 1972                     & 0,07                    & 0\%                                                                    & 8                         & 0,07                      & 0\%                                                                    & 0,06                   \\
                    pmed38                     & 900                      & 400                    & \textbf{1305}            & 1305                     & 1305                     & 0,07                    & 0\%                                                                    & 10                        & 0,07                      & 0\%                                                                    & 0,06                   \\
                    pmed38                     & 900                      & 500                    & \textbf{836}             & 836                      & 836                      & 0,05                    & 0\%                                                                    & 8                         & 0,08                      & 0\%                                                                    & 0,05                   \\\midrule
                    pmed39                     & 900                      & 5                      & \textbf{11069}           & 10938                    & 11069                    & 0,29                    & 0\%                                                                    & 22                        & 3,66                      & 0\%                                                                    & 8,61                   \\
                    pmed39                     & 900                      & 10                     & \textbf{9423}            & 9365                     & 9423                     & 0,31                     & 0\%                                                                    & 18                        & 6,27                      & 0\%                                                                    & 8,01                   \\
                    pmed39                     & 900                      & 20                     & \textbf{7894}            & 7894                     & 7894                     & 0,34                    & 0\%                                                                    & 9                         & 0,52                      & 0\%                                                                    & 0,72                   \\
                    pmed39                     & 900                      & 50                     & \textbf{5941}            & 5937                     & 5941                     & 0,33                    & 0\%                                                                    & 12                        & 0,95                      & 0\%                                                                    & 1,21                    \\
                    pmed39                     & 900                      & 100                    & \textbf{4461}            & 4461                     & 4462                     & 0,24                    & 0\%                                                                    & 9                         & 0,33                      & 0\%                                                                    & 0,36                   \\
                    pmed39                     & 900                      & 200                    & \textbf{2918}            & 2918                     & 2918                     & 0,09                    & 0\%                                                                    & 7                         & 0,14                      & 0\%                                                                    & 0,08                   \\
                    pmed39                     & 900                      & 300                    & \textbf{1968}            & 1968                     & 1968                     & 0,07                    & 0\%                                                                    & 7                         & 0,11                      & 0\%                                                                    & 0,06                   \\
                    pmed39                     & 900                      & 400                    & \textbf{1303}            & 1303                     & 1303                     & 0,08                    & 0\%                                                                    & 10                        & 0,08                      & 0\%                                                                    & 0,06                   \\
                    pmed39                     & 900                      & 500                    & \textbf{821}             & 821                      & 821                      & 0,05                    & 0\%                                                                    & 8                         & 0,08                      & 0\%                                                                    & 0,05                   \\\midrule
                    pmed40                     & 900                      & 5                      & \textbf{12305}           & 12246                    & 12305                    & 0,27                    & 0\%                                                                    & 11                        & 1,41                      & 0\%                                                                    & 3,87                   \\
                    pmed40                     & 900                      & 10                     & \textbf{10491}           & 10439                    & 10491                    & 0,33                    & 0\%                                                                    & 27                        & 2,77                      & 0\%                                                                    & 5,54                   \\
                    pmed40                     & 900                      & 20                     & \textbf{8717}            & 8711                     & 8717                     & 0,39                    & 0\%                                                                    & 13                        & 1,27                      & 0\%                                                                    & 2,39                   \\
                    pmed40                     & 900                      & 50                     & \textbf{6518}            & 6505                     & 6518                     & 0,33                    & 0\%                                                                    & 11                        & 2,94                      & 0\%                                                                    & 4,15                   \\
                    pmed40                     & 900                      & 90                     & \textbf{5128}            & 5128                     & 5128                     & 0,23                    & 0\%                                                                    & 8                         & 0,23                      & 0\%                                                                    & 0,24                   \\
                    pmed40                     & 900                      & 200                    & \textbf{3132}            & 3132                     & 3132                     & 0,11                    & 0\%                                                                    & 8                         & 0,11                      & 0\%                                                                    & 0,09                   \\
                    pmed40                     & 900                      & 300                    & \textbf{2106}            & 2106                     & 2106                     & 0,09                    & 0\%                                                                    & 11                        & 0,09                      & 0\%                                                                    & 0,07                   \\
                    pmed40                     & 900                      & 400                    & \textbf{1398}            & 1398                     & 1398                     & 0,06                    & 0\%                                                                    & 9                         & 0,06                      & 0\%                                                                    & 0,06                   \\
                    pmed40                     & 900                      & 500                    & \textbf{900}             & 900                      & 900                      & 0,04                    & 0\%                                                                    & 7                         & 0,07                      & 0\%                                                                    & 0,05                   \\ \midrule
                    \multicolumn{7}{r|}{Average total time}                                                                                                                                & \multicolumn{3}{|r|}{1,67}                                                                              & \multicolumn{2}{r}{2,69}                                                               \\\bottomrule
                    \end{tabular}%
                    }
                    \caption{Results on ORLIB instances for our method and \Zebra\;on our machine}
                    \label{tab:ORLIB}
                \end{table}

                \begin{table}[H]
                \centering
                \resizebox{0.67\textwidth}{!}{%
                 \begin{tabular}{@{}cccr|rrr|rcr|rr@{}}
                    \toprule
                    \multicolumn{4}{c|}{INSTANCE}                                                                   & \multicolumn{3}{c|}{PHASE 1}                                         & \multicolumn{3}{c|}{PHASE 1 + 2}                                                                                      & \multicolumn{2}{c}{Zebra}                                                             \\ \midrule
                    \multicolumn{1}{c}{Name} & \multicolumn{1}{c}{$N=M$} & \multicolumn{1}{c}{$p$} & \multicolumn{1}{c|}{\begin{tabular}[c]{@{}c@{}}$OPT /$\\ $BKN$\end{tabular}} & \multicolumn{1}{c}{$LB^1$} & \multicolumn{1}{c}{$UB^1$} & \multicolumn{1}{c|}{$T^1$} & \multicolumn{1}{c}{$gap$} & \multicolumn{1}{c}{$iter$} & \multicolumn{1}{c|}{$T^{tot}$} & \multicolumn{1}{c}{$gap$} & \multicolumn{1}{c}{T} \\ \midrule
                rl1304                  & 1304                     & 5                      & \textbf{3099073}         & 3099073                  & 3099073                  & 2,64                    & 0\%                                                                    & 9                         & 2,6                       & 0\%                                                                    & 1232,8                 \\
                rl1304                 & 1304                     & 10                     & \textbf{2134295}         & 2131788                  & 2134295                  & 2,92                    & 0\%                                                                    & 12                        & 15,5                      & 0\%                                                                    & 1060,1                 \\
                rl1304                 & 1304                     & 20                     & \textbf{1412108}         & 1412108                  & 1412108                  & 2,30                    & 0\%                                                                    & 8                         & 2,3                       & 0\%                                                                    & 60,8                   \\
                rl1304                 & 1304                     & 50                     & \textbf{795012}          & 795012                   & 795012                   & 1,48                    & 0\%                                                                    & 9                         & 1,5                       & 0\%                                                                    & 8,3                    \\
                rl1304                & 1304                     & 100                    & \textbf{491639}          & 491507                   & 491788                   & 0,92                    & 0\%                                                                    & 19                        & 2,4                       & 0\%                                                                    & 3,6                    \\
                rl1304                & 1304                     & 200                    & \textbf{268573}          & 268573                   & 268573                   & 0,34                    & 0\%                                                                    & 11                        & 0,5                       & 0\%                                                                    & 0,9                    \\
                rl1304                & 1304                     & 300                    & \textbf{177326}          & 177318                   & 177339                   & 0,28                    & 0\%                                                                    & 12                        & 0,4                       & 0\%                                                                    & 0,5                    \\
                rl1304                & 1304                     & 400                    & \textbf{128332}          & 128332                   & 128332                   & 0,21                    & 0\%                                                                    & 10                        & 0,2                       & 0\%                                                                    & 0,1                    \\
                rl1304                & 1304                     & 500                    & \textbf{97024}           & 97018                    & 97034                    & 0,23                    & 0\%                                                                    & 14                        & 0,3                       & 0\%                                                                    & 0,2                    \\ \midrule
                fl1400                  & 1400                     & 5                      & \textbf{174877}          & 174877                   & 174877                   & 0,73                    & 0\%                                                                    & 7                         & 0,7                       & 0\%                                                                    & 245,3                  \\
                fl1400                 & 1400                     & 10                     & \textbf{100601}          & 100601                   & 100601                   & 0,40                    & 0\%                                                                    & 6                         & 0,4                       & 0\%                                                                    & 71,7                   \\
                fl1400                 & 1400                     & 20                     & \textbf{57191}           & 57191                    & 57191                    & 0,38                    & 0\%                                                                    & 8                         & 0,4                       & 0\%                                                                    & 10,4                   \\
                fl1400                 & 1400                     & 50                     & \textbf{28486}           & 28486                    & 28486                    & 0,36                    & 0\%                                                                    & 8                         & 0,4                       & 0\%                                                                    & 2,5                    \\
                fl1400                & 1400                     & 100                    & \textbf{15962}           & 15961                    & 15962                    & 0,82                    & 0\%                                                                    & 12                        & 2,1                       & 0\%                                                                    & 5,0                    \\
                fl1400                & 1400                     & 200                    & \textbf{8806}            & 8793                     & 8815                     & 0,66                    & 0\%                                                                    & 20                        & 26,8                      & 0\%                                                                    & 304,8                  \\
                fl1400                & 1400                     & 300                    & {\ul \textbf{6109}}      & 6092                     & 6157                     & 0,76                    & 0\%                                                                    & 21                        & 399,7                     & 8,6\%                                                                 & TL                     \\
                fl1400                & 1400                     & 400                    & \textbf{4648}            & 4636                     & 4659                     & 0,45                    & 0\%                                                                    & 46                        & 11964                   & 7,9\%                                                                  & TL                     \\
                fl1400                & 1400                     & 500                    & {3764}            & 3756                     & 3773                     & 0,54                    & 0,08\%                                                                 & 35                        & TL                        & 7,7\%                                                                  & TL                     \\ \midrule
                u1432                   & 1432                     & 5                      & \textbf{1210126}         & 1210126                  & 1210126                  & 1,57                    & 0\%                                                                    & 7                         & 1,6                       & 0\%                                                                    & 323,6                  \\
                u1432                  & 1432                     & 10                     & \textbf{849759}          & 849759                   & 849759                   & 3,50                    & 0\%                                                                    & 7                         & 3,5                       & 0\%                                                                    & 71,2                   \\
                u1432                  & 1432                     & 20                     & \textbf{588766}          & 588720                   & 588767                   & 3,91                    & 0\%                                                                    & 12                        & 5,8                       & 0\%                                                                    & 18,1                   \\
                u1432                  & 1432                     & 50                     & \textbf{362072}          & 361724                   & 362072                   & 4,34                    & 0\%                                                                    & 25                        & 161,8                     & 0\%                                                                    & 128,1                  \\
                u1432                 & 1432                     & 100                    & \textbf{243793}          & 243758                   & 243850                   & 1,89                    & 0\%                                                                    & 12                        & 3,1                       & 0\%                                                                    & 7,3                    \\
                u1432                 & 1432                     & 200                    & \textbf{159887}          & 159867                   & 160084                   & 0,73                    & 0\%                                                                    & 13                        & 2,0                       & 0\%                                                                    & 1,9                    \\
                u1432                 & 1432                     & 300                    & \textbf{123689}          & 123674                   & 123876                   & 0,61                    & 0\%                                                                    & 15                        & 1,0                       & 0\%                                                                    & 2,2                    \\
                u1432                 & 1432                     & 400                    & 103979          & 103411                   & 104102                   & 1,22                    & 0,09\%                                                                 & 32                        & TL                        & 0,5\%                                                                  & TL                     \\
                u1432                 & 1432                     & 500                    & \textbf{93200}           & 93200                    & 93200                    & 0,17                    & 0\%                                                                    & 8                         & 0,3                       & 0\%                                                                    & 0,1                    \\ \midrule
                vm1748                  & 1748                     & 5                      & \textbf{4479421}         & 4479421                  & 4479421                  & 2,23                    & 0\%                                                                    & 7                         & 2,2                       & 0\%                                                                    & 4955,2                 \\
                vm1748                 & 1748                     & 10                     & \textbf{2983645}         & 2983048                  & 2983645                  & 4,64                    & 0\%                                                                    & 14                        & 11,0                      & 0\%                                                                    & 1364,0                 \\
                vm1748                 & 1748                     & 20                     & \textbf{1899680}         & 1899588                  & 1899681                  & 4,76                    & 0\%                                                                    & 15                        & 9,5                       & 0\%                                                                    & 308,6                  \\
                vm1748                 & 1748                     & 50                     & \textbf{1004331}         & 1004325                  & 1004339                  & 2,25                    & 0\%                                                                    & 12                        & 2,8                       & 0\%                                                                    & 21,4                   \\
                vm1748                & 1748                     & 100                    & \textbf{636515}          & 636418                   & 636541                   & 1,60                    & 0\%                                                                    & 15                        & 3,2                       & 0\%                                                                    & 12,6                   \\
                vm1748                & 1748                     & 200                    & \textbf{390350}          & 390350                   & 390350                   & 0,77                    & 0\%                                                                    & 11                        & 0,8                       & 0\%                                                                    & 1,6                    \\
                vm1748                & 1748                     & 300                    & \textbf{286039}          & 286037                   & 286080                   & 0,60                    & 0\%                                                                    & 13                        & 1,0                       & 0\%                                                                    & 1,0                    \\
                vm1748                & 1748                     & 400                    & \textbf{221526}          & 221523                   & 221545                   & 0,48                    & 0\%                                                                    & 13                        & 0,8                       & 0\%                                                                    & 1,0                    \\
                vm1748                & 1748                     & 500                    & \textbf{176986}          & 176977                   & 177103                   & 0,46                    & 0\%                                                                    & 14                        & 0,7                       & 0\%                                                                    & 0,5                    \\ \midrule
                         \multicolumn{7}{r|}{Average total time}                                                                                                                                & \multicolumn{3}{|r|}{8,4}                                                                              & \multicolumn{2}{r}{319,5}                                                               \\ \bottomrule
                \end{tabular}%
                }
                     \caption{Results on medium TSP instances for our method and \Zebra\;on our machine. TL=36000 seconds. The average total time is calculated with the instances in which both methods solve  the instances to optimality. The best computed lower bounds for fl1400(p=500) and u1432(p=400) are 3761 and 103879 respectively}
                \label{tab:TSPmedium}
                \end{table}
 
                \begin{table}[H]
                    \centering
                    \resizebox{0.63\textwidth}{!}{%
                    \begin{tabular}{@{}cccr|rrr|rcr|rr@{}}
                    \toprule
                    \multicolumn{4}{c|}{INSTANCE}                                                                                                         & \multicolumn{3}{c}{PHASE 1}                                                                     & \multicolumn{3}{|c|}{PHASE 1 + 2}                                                                                                                  & \multicolumn{2}{|c}{Zebra}                                                                                \\ \midrule
                    \multicolumn{1}{c}{Name} & \multicolumn{1}{c}{$N=M$} & \multicolumn{1}{c}{p} & \multicolumn{1}{c|}{\begin{tabular}[c]{@{}c@{}}$OPT$\\ BKN\end{tabular}} & \multicolumn{1}{c}{$LB^1$} & \multicolumn{1}{c}{$UB^1$} & \multicolumn{1}{c|}{$T^1$} & \multicolumn{1}{c}{$gap$} & \multicolumn{1}{c}{$iter$} & \multicolumn{1}{c|}{$T^{tot}$} & \multicolumn{1}{c}{$gap$} & \multicolumn{1}{c}{$T$} \\ \midrule
                    d2103                            & 2103                              & 5                               & \textbf{1005136}                  & 1005136                           & 1005136                           & 6                                & 0\%                                                                             & 9                                  & 6                                  & 0\%                                                                             & 4268                            \\
                    d2103                           & 2103                              & 10                              & \textbf{687321}                   & 687264                            & 687321                            & 13                               & 0\%                                                                             & 10                                 & 21                                 & 0\%                                                                             & 1085                            \\
                    d2103                           & 2103                              & 20                              & \textbf{482926}                   & 482798                            & 482926                            & 17                               & 0\%                                                                             & 15                                 & 39                                 & 0\%                                                                             & 448                             \\
                    d2103                           & 2103                              & 50                              & {\underline{302221}}             & 301591                            & 303178                            & 31                               & 0,05\%                                                                          & 30                                 & TL             & 0,4\%                                                                           & TL                            \\
                    d2103                          & 2103                              & 100                             & \textbf{194664}                   & 194407                            & 194994                            & 17                               & 0\%                                                                             & 39                                 & 8033                              & 0\%                                                                             & 16022                           \\
                    d2103                          & 2103                              & 200                             & \textbf{117753}                   & 117735                            & 117778                            & 4                                & 0\%                                                                             & 15                                 & 8                                  & 0\%                                                                             & 5                               \\
                    d2103                          & 2103                              & 300                             & \textbf{90471}                    & 90423                             & 90510                             & 4                                & 0\%                                                                             & 29                                 & 8                                  & 0\%                                                                             & 29                              \\
                    d2103                          & 2103                              & 400                             & \textbf{75324}                    & 75291                             & 75425                             & 2                                & 0\%                                                                             & 20                                 & 6                                  & 0\%                                                                             & 6209                            \\
                    d2103                          & 2103                              & 500                             & \textbf{64006}                    & 63952                             & 64315                             & 2                                & 0\%                                                                             & 31                                 & 12                                 & 0\%                                                                             & 2568                            \\ \midrule
                    pcb3038                          & 3038                              & 5                               & \textbf{1777835}                  & 1777665                           & 1777835                           & 43                               & 0\%                                                                             & 14                                 & 84                                 & 0,1\%                                                                           & TL                           \\
                    pcb3038                         & 3038                              & 10                              & \textbf{1211704}                  & 1211704                           & 1211704                           & 32                               & 0\%                                                                             & 8                                  & 32                                 & 0\%                                                                             & 19526                           \\
                    pcb3038                         & 3038                              & 20                              & \textbf{839494}                   & 839232                            & 839499                            & 86                               & 0\%                                                                             & 19                                 & 235                                & 0,2\%                                                                           & 7329                            \\
                    pcb3038                         & 3038                              & 50                              & \textbf{506339}                   & 506204                            & 506339                            & 44                               & 0\%                                                                             & 13                                 & 136                                & 0\%                                                                             & 1134                            \\
                    pcb3038                        & 3038                              & 100                             & \textbf{351500}                   & 351404                            & 351648                            & 42                               & 0\%                                                                             & 23                                 & 156                                & 0\%                                                                             & 346                             \\
                    pcb3038                        & 3038                              & 150                             & \textbf{280128}                   & 280057                            & 280423                            & 28                               & 0\%                                                                             & 18                                 & 241                                & 0\%                                                                             & 148                             \\
                    pcb3038                        & 3038                              & 200                             & {\ul \textbf{237399}}             & 237328                            & 237578                            & 19                               & 0\%                                                                             & 14                                 & 130                                & 0\%                                                                             & 130                             \\
                    pcb3038                        & 3038                              & 300                             & \textbf{186833}                   & 186793                            & 186906                            & 13                               & 0\%                                                                             & 19                                 & 26                                 & 0\%                                                                             & 60                              \\
                    pcb3038                        & 3038                              & 400                             & \textbf{156276}                   & 156268                            & 156307                            & 6                                & 0\%                                                                             & 10                                 & 11                                 & 0\%                                                                             & 22                              \\
                    pcb3038                        & 3038                              & 500                             & \textbf{134798}                   & 134774                            & 134866                            & 3                                & 0\%                                                                             & 14                                 & 5                                  & 0\%                                                                             & 17                              \\ \midrule
                    fl3795                           & 3795                              & 5                               & \textbf{1052627}                  & 1052627                           & 1052627                           & 12                               & 0\%                                                                             & 7                                  & 13                                 & $\infty$                                                                           & $\blacklozenge$                           \\
                    fl3795                          & 3795                              & 10                              & \textbf{520940}                   & 520940                            & 520940                            & 7                                & 0\%                                                                             & 8                                  & 7                                  & 0\%                                                                             & 4410                            \\
                    fl3795                          & 3795                              & 20                              & \textbf{319722}                   & 319722                            & 319722                            & 7                                & 0\%                                                                             & 10                                 & 7                                  & 0\%                                                                             & 2671                            \\
                    fl3795                          & 3795                              & 50                              & \textbf{150940}                   & 150940                            & 150940                            & 5                                & 0\%                                                                             & 12                                 & 5                                  & 0\%                                                                             & 193                             \\
                    fl3795                         & 3795                              & 100                             & \textbf{88299}                    & 88299                             & 88299                             & 6                                & 0\%                                                                             & 11                                 & 6                                  & 0\%                                                                             & 45                              \\
                    fl3795                         & 3795                              & 150                             & \textbf{65868}                    & 65840                             & 65904                             & 14                               & 0\%                                                                             & 52                                 & 217                                & 0\%                                                                             & 1825                            \\
                    fl3795                         & 3795                              & 200                             & \textbf{53928}                    & 53913                             & 54013                             & 12                               & 0\%                                                                             & 66                                 & 2732                               & 0\%                                                                             & 2501                            \\
                    fl3795                         & 3795                              & 300                             & \textbf{39586}                    & 39578                             & 39661                             & 7                                & 0\%                                                                             & 35                                 & 2772                               & 0\%                                                                             & 3061                            \\
                    fl3795                         & 3795                              & 400                             & \textbf{31354}                    & 31348                             & 31472                             & 6                                & 0\%                                                                             & 36                                 & 454                                & 0\%                                                                             & 527                             \\
                    fl3795                         & 3795                              & 500                             & \textbf{25976}                    & 25976                             & 25988                             & 6                                & 0\%                                                                             & 16                                 & 6                                  & 0\%                                                                             & 2                               \\ \midrule
                    rl5934                          & 5934                              & 10                              & \textbf{\underline{9792218}}            & 9786688                           & 9792218                           & 545                              & 0\%                                                                             & 17                                 & 3162                               &    $\infty$                                                                             &   $\blacklozenge$                               \\
                    rl5934                          & 5934                              & 20                              & \textbf{\underline{6716215}}            & 6713214                           & 6716228                           & 535                              & 0\%                                                                             & 30                                 & 19337                              &    $\infty$                                                                             &     $\blacklozenge$                               \\
                    rl5934                          & 5934                              & 50                              & \underline{4030007}                  & 4026935                           & 4032425                           & 464                              & 0,03\%                                                                          & 32                                 & {TL}             & 0,1\%                                                                           & TL                              \\
                    rl5934                         & 5934                              & 200                             & \textbf{1805530}                & 1805029                           & 1807763                         & 104                              & 0\%                                                                             & 25                                 & 1479                               & 0\%                                                                             & 3816                            \\
                    rl5934                         & 5934                              & 300                             & \textbf{1392419}                & 1392304                           & 1392709                         & 89                               & 0\%                                                                             & 15                                 & 89                                 & 0\%                                                                             & 235                             \\
                    rl5934                         & 5934                              & 400                             & \textbf{1143940}                & 1143649                           & 1145342                         & 575                              & 0\%                                                                             & 30                                 & 575                                & 0\%                                                                             & 1110                            \\
                    rl5934                         & 5934                              & 500                             & \textbf{972799}                  & 972741                            & 973712                           & 57                               & 0\%                                                                             & 16                                 & 57                                 & 0\%                                                                             & 70                              \\ 
                    rl5934                         & 5934                              & 600                             & \textbf{847301}                  & 847232                            & 847769                           & 34                               & 0\%                                                                             & 19                                 & 34                                 & 0\%                                                                             & 77                              \\
                    rl5934                         & 5934                              & 700                             & \textbf{751131}                  & 751054                            & 751569                           & 23                               & 0\%                                                                             & 15                                 & 23                                 & 0\%                                                                             & 88                              \\
                    rl5934                         & 5934                              & 800                             & \textbf{675958}                  & 675884                            & 676248                           & 26                               & 0\%                                                                             & 15                                 & 26                                 & 0\%                                                                             & 58                              \\
                    rl5934                         & 5934                              & 900                             & \textbf{612629}                  & 612574                            & 612879                           & 18                               & 0\%                                                                             & 14                                 & 18                                 & 0\%                                                                             & 33                              \\
                    rl5934                        & 5934                              & 1000                            & \textbf{558167}                  & 558087                            & 558311                           & 85                               & 0\%                                                                             & 37                                 & 85                                 & 0\%                                                                             & 731                             \\
                    rl5934                        & 5934                              & 1100                            & \textbf{511192}                  & 511138                            & 511453                           & 25                               & 0\%                                                                             & 22                                 & 25                                 & 0\%                                                                             & 20                              \\
                    rl5934                        & 5934                              & 1200                            & \textbf{469747}                  & 469711                            & 469943                           & 16                               & 0\%                                                                             & 19                                 & 16                                 & 0\%                                                                             & 11                              \\
                    rl5934                        & 5934                              & 1300                            & \textbf{433060}                  & 433014                            & 433300                           & 17                               & 0\%                                                                             & 22                                 & 17                                 & 0\%                                                                             & 27                              \\
                    rl5934                       & 5934                              & 1400                            & \textbf{401370}                  & 401356                            & 401542                           & 14                               & 0\%                                                                             & 16                                 & 14                                 & 0\%                                                                             & 6                               \\
                    rl5934                        & 5934                              & 1500                            & \textbf{373566}                  & 373566                            & 373566                           & 10                               & 0\%                                                                             & 18                                 & 10                                 & 0\%                                                                             & 2 \\ \midrule
                         \multicolumn{7}{r|}{Average total time}                                                                                                                                & \multicolumn{3}{|r|}{470}                                                                              & \multicolumn{2}{r}{2022}                                                               \\ \bottomrule
                \end{tabular}
                }
                 \caption{Results on large TSP instances for our method and \Zebra\;on our machine. TL=36000 seconds. The average total time is calculated with the instances in which both methods solve the instances to optimality. 
                 }
                    \label{tab:TSPlarge}
                \end{table}
                
                \begin{table}[H]
                    \centering
                    \resizebox{0.69\textwidth}{!}{%
                    \begin{tabular}{@{}cccr|rrr|rcr|rr@{}}
                    \toprule
                    \multicolumn{4}{c|}{INSTANCE}                                                                                                         & \multicolumn{3}{c}{PHASE 1}                                                                     & \multicolumn{3}{|c|}{PHASE 1 + 2}                                                                                                                  & \multicolumn{2}{|c}{Zebra}                                                                                \\ \midrule
                    \multicolumn{1}{c}{Name} & \multicolumn{1}{c}{$N=M$} & \multicolumn{1}{c}{p} & \multicolumn{1}{c|}{\begin{tabular}[c]{@{}c@{}}$OPT /$\\ $BKN$\end{tabular}} & \multicolumn{1}{c}{$LB^1$} & \multicolumn{1}{c}{$UB^1$} & \multicolumn{1}{c|}{$T^1$} & \multicolumn{1}{c}{$gap$} & \multicolumn{1}{c}{$iter$} & \multicolumn{1}{c|}{$T^{tot}$} & \multicolumn{1}{c}{$gap$} & \multicolumn{1}{c}{$T$} \\ \midrule
                    usa13509              & 13509                    & 300                    & \textbf{59340915}        & 59334913                 & 59346783                 & 506                     & 0\%                                                                    & 14                        & 1663                      & 0\%                                                                    & 18760                  \\
                    usa13509              & 13509                    & 400                    & \textbf{50538905}        & 50533013                 & 50575463                 & 333                     & 0\%                                                                    & 16                        & 2578                      & 0\%                                                                    & 23677                  \\
                    usa13509              & 13509                    & 500                    & \textbf{44469860}        & 44463038                 & 44499566                 & 289                     & 0\%                                                                    & 20                        & 3715                      & 0\%                                                                    & 25105                  \\
                    usa13509              & 13509                    & 600                    & {\ul \textbf{39952138}}  & 39944049                 & 39991088                 & 301                     & 0\%                                                                & 21                        & 16033                     & $\infty$                                                                     & $\blacklozenge$                  \\
                    usa13509              & 13509                    & 700                    & {\ul \textbf{36469603}}  & 36463603                 & 36512930                 & 203                     & 0\%                                                                    & 22                        & 1953                      & $\infty$                                                                     & $\blacklozenge$                  \\
                    usa13509              & 13509                    & 800                    & \textbf{33635127}        & 33631192                 & 33672848                 & 215                     & 0\%                                                                    & 27                        & 1012                      & 0\%                                                                    & 6007                   \\
                    usa13509              & 13509                    & 900                    & {\ul \textbf{31275114}}  & 31269089                 & 31299760                 & 177                     & 0\%                                                                    & 18                        & 5043                      & $\infty$                                                                     & $\blacklozenge$                  \\
                    usa13509             & 13509                    & 1000                   & \textbf{29268216}        & 29262339                 & 29309009                 & 151                     & 0\%                                                                    & 20                        & 823                       & $\infty$                                                                     & $\blacklozenge$                  \\
                    usa13509             & 13509                    & 2000                   & \textbf{18230856}        & 18229432                 & 18238229                 & 38                      & 0\%                                                                    & 18                        & 115                       & 0\%                                                                    & 584                    \\
                    usa13509             & 13509                    & 3000                   & \textbf{13098935}        & 13097929                 & 13101469                 & 25                      & 0\%                                                                    & 23                        & 46                        & 0\%                                                                    & 1674                   \\
                    usa13509             & 13509                    & 4000                   & \textbf{9905715}         & 9905071                  & 9910848                  & 18                     & 0\%                                                                    & 16                        & 27                        & 0\%                                                                    & 166                    \\
                    usa13509             & 13509                    & 5000                   & \textbf{7608605}         & 7608242                  & 7611958                  & 16                      & 0\%                                                                    & 23                        & 27                        & 0\%                                                                    & 86                     \\ \midrule
                    sw24978              & 24978                    & 1000                   & {\ul \textbf{1841723}}   & 1841613                  & 1844801                  & 515                     & 0\%                                                                    & 21                        & 4211                      & 0\%                                                                    & 30796                  \\
                    sw24978              & 24978                    & 2000                   & {\ul \textbf{1197278}}   & 1197231                  & 1198464                  & 174                     & 0\%                                                                    & 16                        & 332                       & $\infty$                                                                    & $\blacklozenge$                  \\
                    sw24978              & 24978                    & 3000                   & {\ul \textbf{911361}}    & 911308                   & 911988                   & 124                     & 0\%                                                                    & 13                        & 196                       & 0\%                                                                    & 3614                   \\
                    sw24978              & 24978                    & 4000                   & {\ul \textbf{737645}}    & 737602                   & 738045                   & 79                     & 0\%                                                                    & 15                        & 125                       & $\infty$                                                                     & $\blacklozenge$                  \\
                    sw24978              & 24978                    & 5000                   & {\ul \textbf{617637}}    & 617593                   & 618096                   & 66                     & 0\%                                                                    & 18                        & 101                      & $\infty$                                                                     & $\blacklozenge$                  \\
                    sw24978              & 24978                    & 6000                   & {\ul \textbf{527336}}    & 527307                   & 527716                   & 60                     & 0\%                                                                    & 17                        & 92                       & 0\%                                                                    & 5188                   \\
                    sw24978              & 24978                    & 7000                   & {\ul \textbf{455716}}    & 455696                   & 456074                   & 53                      & 0\%                                                                    & 16                        & 85                       & 0\%                                                                    & 3973                   \\
                    sw24978              & 24978                    & 8000                   & {\ul \textbf{397217}}    & 397153                   & 397540                   & 38                      & 0\%                                                                    & 18                        & 68                       & $\infty$                                                                    & $\blacklozenge$                  \\
                    sw24978              & 24978                    & 9000                   & {\ul \textbf{347376}}    & 347322                   & 347621                   & 37                      & 0\%                                                                    & 22                        & 100                      & $\infty$                                                                     & $\blacklozenge$                  \\
                    sw24978             & 24978                    & 10000                  & {\ul \textbf{305998}}    & 305932                   & 306094                   & 27                      & 0\%                                                                    & 20                        & 49                        & $\infty$                                                                     & $\blacklozenge$                  \\ \midrule
                    ch71009             & 71009                    & 10000                  & {\ul {4274688}}         & 4273680                        & 4417777                        & 5082                       & 0,007\%                                                                    & 32                         & TL                         & $\infty$                                                                       & $\blacklozenge$                       \\
                    ch71009             & 71009                    & 20000                  & {\ul \textbf{2377760}}   & 2377409                  & 2419539                  & 687                     & 0\%                                                                    & 40                        & 3560                      & $\infty$                                                                        & $\blacklozenge$                       \\
                    ch71009             & 71009                    & 30000                  & {\ul \textbf{1464151}}   & 1464015                  & 1473517                  & 438                    & 0\%                                                                    & 27                        & 814                      & $\infty$                                                                       & $\blacklozenge$                       \\
                    ch71009             & 71009                    & 40000                  & {\ul \textbf{879336}}    & 879272                   & 881997                   & 221                     & 0\%                                                                    & 17                        & 445                       &  $\infty$                                                                      & $\blacklozenge$                       \\
                    ch71009             & 71009                    & 50000                  & {\ul \textbf{463553}}    & 463544                   & 463904                   & 135                     & 0\%                                                                    & 24                        &  231                      & 0\%                                                                    & 653                   \\
                    ch71009             & 71009                    & 60000                  & {\ul \textbf{167565}}    & 167558                   & 167789                   & 50                     & 0\%                                                                    & 31                        & 102                       & 0\%                                                                    & 331                   \\ \midrule
                    pla85900            & 85900                    & 10000                  & {\ul {166855420}} & 166627292                & 187017656                & 2619                    & 0,124\%                                                               & 30                        & TL    &   $\infty$                                                                     & $\blacklozenge$                       \\
                    pla85900            & 85900                    & 20000                  & {\ul {108208172}} & 107246411                & 121125450                & 1257                    & 0,829\%                                                               & 20                        & TL    &   $\infty$                                                                     & $\blacklozenge$                       \\
                    pla85900            & 85900                    & 30000                  & {\ul {86945099}}  & 86944715                 & 877780783                 & 371                    & 0,0004\%                                                               & 70                        & TL    &   $\infty$                                                                     & $\blacklozenge$                       \\
                    pla85900            & 85900                    & 40000                  & {\ul {69944760}}  & 69944715                 & 69979648                 & 143                     & 0,0001\%                                                               & 60                        & TL    &   $\infty$                                                                     & $\blacklozenge$                       \\
                    pla85900            & 85900                    & 50000                  & {\ul \textbf{52944715}}  & 52944715                 & 52945981                 & 47                     & 0\%                                                               & 38                        & 32357   &  $\infty$                                                                      & $\blacklozenge$                       \\
                    pla85900            & 85900                    & 60000                  & {\ul \textbf{35944715}}  & 35944715                 & 35945264                 & 67                     & 0\%                                                                    & 20                        & 510                      &      $\infty$                                                                  & $\blacklozenge$                       \\
                    pla85900            & 85900                    & 70000                  & {\ul \textbf{18977475}}  & 18977475                 & 18977475                 & 76                     & 0\%                                                                    & 13                        & 76                       & 0\%                                                                    & 122                    \\
                    pla85900            & 85900                    & 80000                  & {\ul \textbf{4512752}}   & 4512752                  & 4512752                  & 13                     & 0\%                                                                    & 20                        & 13                       & 0\%                                                                    & 97                    \\ \midrule
                    usa115475           & 115475                   & 20000                  & {\ul 5287409}            & 5286659                  & 5385388                  & 3736                    & 0,002\%                                                               & 37                        & TL                     &   $\infty$                                                                     & $\blacklozenge$                       \\
                    usa115475           & 115475                   & 30000                  & {\ul \textbf{3815620}}   & 3815143                  & 3848533                  & 1317                    & 0\%                                                                    & 34                        & 9529                      &   $\infty$                                                                     & $\blacklozenge$                       \\
                    usa115475           & 115475                   & 40000                  & {\ul \textbf{2876909}}   & 2876603                  & 2904492                  & 1363                    & 0\%                                                                    & 32                        & 4412                      &  $\infty$                                                                      & $\blacklozenge$                       \\
                    usa115475           & 115475                   & 50000                  & {\ul \textbf{2189144}}   & 2188903                  & 2200969                  & 1107                    & 0\%                                                                    & 28                        & 3097                      &  $\infty$                                                                      & $\blacklozenge$                       \\
                    usa115475           & 115475                   & 60000                  & {\ul \textbf{1651400}}   & 1651234                  & 1657118                  & 774                    & 0\%                                                                    & 25                        & 1501                      &  $\infty$                                                                      & $\blacklozenge$                       \\
                    usa115475           & 115475                   & 70000                  & {\ul \textbf{1214299}}   & 1214177                  & 1217251                  & 597                    & 0\%                                                                    & 17                        & 964                      &  $\infty$                                                                      & $\blacklozenge$                       \\
                    usa115475           & 115475                   & 80000                  & {\ul \textbf{851481}}    & 851422                   & 852851                   & 434                    & 0\%                                                                    & 24                        & 719                      &  $\infty$                                                                      & $\blacklozenge$                       \\
                    usa115475           & 115475                   & 90000                  & {\ul \textbf{548097}}    & 548076                   & 548560                   & 270                    & 0\%                                                                    & 18                        & 468                      & $\infty$                                                                   & $\blacklozenge$   \\
                    \midrule
                    \multicolumn{7}{r|}{Average total time}                             & \multicolumn{3}{|r|}{887}                                       & \multicolumn{2}{r}{7552}                                                   \\ \bottomrule
                    \end{tabular}%
                    }
                     \caption{Results on huge TSP instances for our method and \Zebra\;on our machine. TL=36000 seconds. The average total time is calculated with the instances in which both methods solve the instances to optimality}
                    \label{tab:TSPhuge}
                    \end{table}
                    
                \begin{table}[H]
                    \centering
                    \resizebox{0.64\textwidth}{!}{%
                    \begin{tabular}{@{}cccr|rrr|rcr|rr@{}}
                    \toprule
                    \multicolumn{4}{c|}{INSTANCE}                                                                                                         & \multicolumn{3}{c}{PHASE 1}                                                                     & \multicolumn{3}{|c|}{PHASE 1 + 2}                                                                                                                  & \multicolumn{2}{|c}{PopStar}                                                                                \\ \midrule
                    \multicolumn{1}{c}{Name} & \multicolumn{1}{c}{$N=M$} & \multicolumn{1}{c}{p} & \multicolumn{1}{c|}{\begin{tabular}[c]{@{}c@{}}$OPT /$\\ $BKN$\end{tabular}} & \multicolumn{1}{c}{$LB^1$} & \multicolumn{1}{c}{$UB^1$} & \multicolumn{1}{c|}{$T^1$} & \multicolumn{1}{c}{$gap$} & \multicolumn{1}{c}{$iter$} & \multicolumn{1}{c|}{$T^{tot}$} & \multicolumn{1}{c}{$UB^h$} & \multicolumn{1}{c}{$T$} \\ \midrule
                    rw100                           & 100                               & 10                              & \textbf{530}                                                                     & 475                              & 530                              & 0,03                                & 0\%                                   & 45                                 & 10,29                               & \textbf{530}                                                                     & 1                                 \\
                    rw100                           & 100                               & 20                              & \textbf{277}                                                                     & 274                              & 277                              & 0,02                                & 0\%                                   & 9                                  & 0,24                                & \textbf{277}                                                                     & 7                                   \\
                    rw100                           & 100                               & 30                              & \textbf{213}                                                                     & 213                              & 213                              & 0,01                                & 0\%                                   & 9                                  & 0,01                                & \textbf{213}                                                                     & 0,5                                 \\
                    rw100                           & 100                               & 40                              & \textbf{187}                                                                     & 187                              & 187                              & 0,01                                & 0\%                                   & 7                                  & 0,01                                & \textbf{187}                                                                     & 0,5                                 \\
                    rw100                           & 100                               & 50                              & \textbf{172}                                                                     & 172                              & 172                              & 0,00                                & 0\%                                   & 7                                  & 0,01                                & \textbf{172}                                                                     & 0,4                                 \\ \midrule
                    rw250                           & 250                               & 10                              & \textbf{3691}                                                                    & 2811                             & 3698                             & 0,35                                & 1,9\%                                 & 92                                 & TL                               & {3698}                                                                    & 10                                \\
                    rw250                           & 250                               & 25                              & \textbf{\underline{1360}}                                                              & 1216                             & 1364                         & 0,23                                & 0\%                                   & 103                                & 20181                               & 1364                                                                             & 6                                 \\
                    rw250                           & 250                               & 50                              & \textbf{713}                                                                     & 699                              & 713                              & 0,17                                & 0\%                                   & 21                                 & 7,81                               & \textbf{713}                                                                     & 4                                 \\
                    rw250                           & 250                               & 75                              & \textbf{523}                                                                     & 523                              & 523                              & 0,04                                & 0\%                                   & 7                                  & 0,06                                & \textbf{523}                                                                     & 3                                 \\
                    rw250                          & 250                               & 100                             & \textbf{444}                                                                     & 444                              & 444                              & 0,02                                & 0\%                                   & 8                                  & 0,02                                & \textbf{444}                                                                     & 2                                \\
                    rw250                          & 250                               & 125                             & \textbf{411}                                                                     & 411                              & 411                              & 0,02                                & 0\%                                   & 5                                  & 0,04                                & \textbf{411}                                                                     & 2                                   \\ \midrule
                    rw500                           & 500                               & 10                              & 16108                                                                            & 11012                            & 16144                            & 2,91                                & 23,2\%                                & 85                                 &     TL                           & 16144                                                                            & 77                                \\
                    rw500                           & 500                               & 25                              & 5683                                                                             & 4403                             & 5716                             & 1,71                                & 17,2\%                                & 72                                 & TL                               & 5716                                                                             & 47                                \\
                    rw500                           & 500                               & 50                              & 2627                                                                             & 2321                             & 2627                             & 1,03                                & 6,6\%                                 & 52                                 & TL                               & 2635                                                                             & 28                                \\
                    rw500                           & 500                               & 75                              & 1757                                                                             & 1672                             & 1757                             & 0,64                                & 1,3\%                                 & 27                                 & TL                               & 1757                                                                             & 21                                \\
                    rw500                          & 500                               & 100                             & \textbf{1379}                                                                    & 1353                             & 1382                             & 0,35                                & 0\%                                   & 46                                 & 2302                                & \textbf{1382}                                                                    & 20                               \\
                    rw500                          & 500                               & 150                             & \textbf{1024}                                                                    & 1024                             & 1024                             & 0,07                                & 0\%                                   & 8                                  & 0,09                                & \textbf{1024}                                                                    & 15                                \\
                    rw500                          & 500                               & 250                             & \textbf{833}                                                                     & 833                              & 833                              & 0,04                                & 0\%                                   & 9                                  & 0,05                                & \textbf{833}                                                                     & 12                                \\ \midrule
                    rw1000                          & 1000                              & 10                              & 68136                                                                            & 44697                            & 68136                            & 23,71                               & 38,1\%                                & 63                                 & TL                               & 68136                                                                            & 256                               \\
                    rw1000                          & 1000                              & 25                              & 24964                                                                            & 17387                            & 25042                            & 12,05                               & 34,1\%                                & 53                                 & TL                               & 25042                                                                            & 294                               \\
                    rw1000                          & 1000                              & 50                              & 11334                                                                            & 8760                             & 11328                            & 9,54                                & 23,7\%                                & 59                                 & TL                               & 11360                                                                            & 169                               \\
                    rw1000                          & 1000                              & 75                              & 7207                                                                             & 5998                             & 7223                             & 7,79                                & 16,6\%                                & 63                                 & TL                               & 7223                                                                            & 160                              \\
                    rw1000                         & 1000                              & 100                             & 5234                                                                             & 4631                             & 5233                             & 5,93                                & 9,9\%                                 & 44                                 & TL                               & 5259                                                                             & 110                               \\
                    rw1000                         & 1000                              & 200                             & 2710                                                                             & 2664                             & 2710                             & 2,12                                & 0,5\%                                 & 25                                 & TL                               & 2710                                                                             & 100                               \\
                    rw1000                         & 1000                              & 300                             & \textbf{2018}                                                                    & 2017                             & 2018                             & 0,28                                & 0\%                                   & 13                                 & 0,41                                & \textbf{2018}                                                                    & 72                                \\
                    rw1000                         & 1000                              & 400                             & \textbf{1734}                                                                    & 1734                             & 1734                             & 0,12                                & 0\%                                   & 9                                  & 0,12                                & \textbf{1734}                                                                    & 74                                \\
                    rw1000                         & 1000                              & 500                             & \textbf{1614}                                                                    & 1614                             & 1614                             & 0,11                                & 0\%                                   & 8                                  & 0,17                                & \textbf{1614}                                                                    & 56                               \\ \bottomrule
                     \end{tabular}
                     }
                    \caption{Results on RW instances for our exact method and \Popstar\;heuristic in our machine. TL=36000 seconds.}
                    \label{tab:RW}
            \end{table}
    
                \begin{table}[H]
                    \centering
                    \resizebox{0.66\textwidth}{!}{%
                    \begin{tabular}{@{}cccr|rrr|rrr|rr@{}}
                    \toprule
                    \multicolumn{4}{c|}{INSTANCE}                                                                                                         & \multicolumn{3}{c}{PHASE 1}                                                                     & \multicolumn{3}{|c|}{PHASE 1 + 2}                                                                                                                  & \multicolumn{2}{|c}{AvellaB\&C}                                                                                \\ \midrule
                    \multicolumn{1}{c}{Name} & \multicolumn{1}{c}{N=M} & \multicolumn{1}{c}{p} & \multicolumn{1}{c|}{\begin{tabular}[c]{@{}c@{}}$OPT /$\\ $BKN$\end{tabular}} & \multicolumn{1}{c}{$LB^1$} & \multicolumn{1}{c}{$UB^1$} & \multicolumn{1}{c|}{$T^1$} & \multicolumn{1}{c}{$gap$} & \multicolumn{1}{c}{$iter$} & \multicolumn{1}{c|}{$T^{tot}$} & \multicolumn{1}{c}{$gap$} & \multicolumn{1}{c}{$T$} \\ \midrule
                    BD3773                         & 3773         & 5     & \textbf{726954998,4}                                                             & 715785543,5                     & 748669824,0                     & 5,83                               & 0\%                              & 10104                             & 39                                 & 0\%                              & 1540                            \\
                    BD3773                         & 3773         & 6     & \textbf{685812258,0}                                                             & 673317265,9                     & 720632505,6                     & 9,92                               & 0\%                              & 12313                             & 131                                & 0\%                              & 41551                           \\
                    BD3773                        & 3773         & 7     & \textbf{651930471,0}                                                             & 636565701,5                     & 727428978,0                     & 10,98                              & 0\%                              & 16012                             & 562                                & 0\%                              & 216851                          \\
                    BD3773                         & 3773         & 8     & {\ul \textbf{620886605,4}}                                                       & 606599816,1                     & 1157543276,4                    & 19,78                              & 0\%                              & 19245                             & 1018                               & 1,7\%                            & 329053                          \\
                    BD3773                         & 3773         & 9     & {\ul \textbf{595955799,0}}                                                       & 581022150,3                     & 649313415,0                     & 18,29                              & 0\%                              & 16956                             & 1462                               & 2,6\%                            & TL2                          \\
                    BD3773                        & 3773         & 10     & {\ul \textbf{574634206,8}}                                                       & 559096162,3                     & 633383307,0                     & 23,29                              & 0\%                              & 19429                             & 3164                               & 2,8\%                            & TL2                          \\
                    BD3773                        & 3773         & 12     & {\ul \textbf{536700087,0}}                                                       & 522614063,7                     & 605415767,4                     & 27,33                              & 0\%                              & 18635                             & 9610                               & 3,1\%                            & TL2                          \\
                    BD3773                        & 3773         & 13     & {\ul \textbf{521375065,2}}                                                       & 507136693,1                     & 581851884,6                     & 29,04                              & 0\%                              & 20439                             & 12848                              & 3,2\%                            & TL2                          \\
                    BD3773                        & 3773         & 14     & {\ul \textbf{507510543,6}}                                                       & 493051932,7                     & 550267457,4                     & 35,42                              & 0\%                              & 20883                             & 33169                              & 3,1\%                            & TL2  \\ \bottomrule                       
                    \end{tabular}
                    }
                     \caption{Results on ODM instances for our method on our machine and the results of \Avella\;reported in \cite{avella2007computational} TL2=360000 seconds.}
                    \label{tab:ODM}
                    \end{table}

    \section{Conclusions }\label{sec:conclusions}

        The $p$-median problem is a well-studied operations research problem in which we have to choose $p$ sites among $M$ to allocate $N$ clients in order to minimize the sum of their allocation distances. This problem has various applications and several heuristics methods have been proposed to solve it. However, its exact resolution remains a challenge for large instances. The most effective approach in the literature  is able to solve a few instances with 85900 clients and sites. 
        
        We performed a Benders decomposition of the most efficient formulation of the $p$-median problem.  The efficiency of our decomposition comes from a polynomial algorithm for the resolution of the sub-problems in conjunction with several implementation improvements in a two-phase resolution. In the first phase, the integrity constraints are relaxed and in the second phase, the problem is solved in an efficient \textit{branch-and-Benders-cut} approach. 
        
        Our approach outperforms other state-of-the-art methods. We solve instances  having up to 115475 clients and sites from the OR and TSP libraries.     We also tested our decomposition on other $p$-median instances: RW instances which do not satisfy triangle inequalities and ODM instances in which there are allocation prohibitions between certain customers and sites. For the RW instances. we were able to solve instances of up to 1000 clients with a large $p$ value. For ODM instances with 3773 clients, we solve all instances within a 10-hours time limit. 
        
        One of the perspectives of this research is to exploit these results on other families of location problems. It is also expected to use other branching strategies that allow a greater efficiency during the development of the \textit{branch-and-Benders-cut} algorithm. In this sense, an implementation in which the distance matrix is not stored in memory and the corresponding distances are calculated only when needed could allow better scalability of our method.
    
        \section*{Acknowledgments}
    
        \noindent The authors would like to thank Sergio Garcia for providing the code used in~\cite{Garcia2011}.
        
        \noindent This work was funded by the National Agency for Research and Development of Chile - ANID (Scholarship Phd. Program 2019-72200492).
    
        \newpage
        \bibliography{bib_p_median_Benders}

\end{document}